\documentclass[a4paper]{article}
\usepackage{amssymb,amsfonts,amsthm}\usepackage[tbtags]{amsmath}

\usepackage{mathrsfs}
\usepackage{tikz}
\usetikzlibrary{calc}
\tikzstyle{wbl}=[line width=1pt,draw=white,double=black,double distance=0.4pt]
\usepackage{graphicx}
\usepackage{tikz-cd}

\usepackage{hyperref}

\usepackage{titlesec}
\titleformat*{\section}{\normalsize\bfseries}
\titleformat*{\subsection}{\normalsize\bfseries}

\renewcommand{\theequation}{\arabic{section}.\arabic{equation}}
\makeatletter
\@addtoreset{equation}{section}
\makeatother
\newcommand{\Angles}[2]{\langle\!\langle#1,#2\rangle\!\rangle}
\newcommand{\angles}[2]{\langle#1,#2\rangle}

\newcommand{\Proof}{\noindent\emph{Proof}}
  \usepackage{authblk}

\title{Infinitely many Brake orbits of Tonelli Hamiltonian systems on the cotangent bundle} 
\author{\small Duanzhi Zhang\thanks{Corresponding author, email: 
zhangdz@nankai.edu.cn}}
\author{\small Zhihao Zhao\thanks{Email: zhaozhihao@mail.nankai.edu.cn}}
\affil{\small School of Mathematical Sciences, Nankai University, China}
\date{}
\newtheorem{thm}{Theorem}[section]
\newtheorem{lem}[thm]{Lemma}
\newtheorem{prop}[thm]{Proposition}
\newtheorem{cor}[thm]{Corollary}
\theoremstyle{definition}
\newtheorem{defe}{Definition}[section]
\theoremstyle{remark}
\newtheorem{rmk}{Remark}[section]
\newtheorem*{clm}{Claim}
\newtheorem*{conditionsofSPL}{Conditions of generalized splitting Lemma}
\begin{document}

% \makeatletter % change default title style
% \renewcommand*\maketitle{%
%     \begin{center}% 居中标题
%         \bfseries % 默认粗体
%         {\Large \@title \par} % LARGE字号
%         \vskip 1em% %%%  标题下面只有1em的缩进或margin
%         {\global\let\author\@empty}%
%         {\global\let\date\@empty}%
%        % \thispagestyle{empty} %  不设置页面样式
%     \end{center}%
%   \setcounter{footnote}{0}%
% }
% \makeatother
\maketitle

\begin{abstract}
  We prove that on the twisted cotangent bundle of a closed manifold with an exact magnetic form,
 a Hamiltonian system 
of a time-dependent Tonelli Hamiltonian possesses
infinitely many brake orbits.
More precisely, by applying Legendre transform we show that there are infinitely many symmetric orbits of the dual Euler-Lagrange system on the configuration space.
This result contains an assertion for the existence of infinitely many symmetric orbits of Tonelli Euler-Lagrange systems given by G. Lu at the end 
of \cite[Remark 6.1]{lu_corrigendum_2009}. In this paper, we will present a complete proof of this assertion. 
% if the anti-symplectic invoution with nonempty fixed points is defined and preserves the Hamilition function. 

% if the anti-symplectic invoution with nonempty fixed points is defined and preserves the Hamilition function. 
%The statement extends a result due to Guangcun Lu for convex quadratic-growth 
% Hamilitions on cotangent bundles with the standard symplectic form. 

\end{abstract}

 \tableofcontents

\newpage
\pagenumbering{arabic}

 \section{Introduction}
Let $M$ be a closed $N$-dimensional smooth manifold.
 We denote by the triple $(T^*M,\omega,R)$
 the symplectic manifold $(T^*M,\omega)$ associated with an anti-symplectic involution $R$, i.e. $R^*\omega=-\omega$ and $R^2=id$.
Let $H: S_\tau\times T^*M\rightarrow \mathbb{R}$ be a smooth $\tau$-periodic Hamiltonian on symplectic manifold
 $(T^*M,\omega)$. In this paper, we focus on $H$ with the following conditions:

(H1) the fiberwise Hessian of $H$ is positive definite, i.e.
\begin{equation*}
  \sum_{i,j=1}^N\frac{\partial^2H}{\partial p_i\partial p_j}(t,q,p)w_iw_j>0
\end{equation*}
 for all $(t,q,p)\in S_\tau\times T^*M$ and nonzero $w=\sum_{i=1}^N w_i\tfrac{\partial}{\partial p_i}\in T_pT_q^*M.$

 (H2) $H$ is fiberwise superlinear, i.e.
\begin{equation*}
  \lim_{|p|_q\rightarrow\infty}\frac{H(t,q,p)}{|p|_q}=\infty,
\end{equation*}
for all $(t,q)\in S_\tau\times M.$

(H3)  $R^*H(-t,\cdot)=H(t,\cdot).$ 

(H1)(H2) are the so-called Tonelli conditions, and we say $H$ a Tonelli Hamiltonian.
Besides, if $H$ satisfies (H3), we say $H$ a Hamiltonian defined on $(T^*M,\omega,R)$.
 A $\tau$-periodic brake orbit
is defined as the solution of the following Hamilition system 
\begin{equation}\label{sec1:02}
  \begin{aligned}
  \dot{x}(t)&=X_H(t,x(t)),\\
  Rx(t)&=x(-t),\\
  x(t)&=x(t+\tau),
\end{aligned}
\end{equation}
where $X_H$ is a Hamiltonian vector field defined by $i_{X_H}\omega=-\text{d}H$.

If $H=K+V$ is the sum of kinetic energy and the potential energy, which is symmetric in the momenta, then $R_0^*H=H$ and
    $H$ defines the classical mechanical
   system on $(T^*M,\omega_0,R_0)$,
   where $\omega_0=\text{d}q\wedge\text{d}p,\ p\in T_q^*M$ and $R_0: (q,p)\rightarrow(q,-p)$ is the standard anti-symplectic involution on $T^*M$ with zero 
   section $M_0$ as the fixed point set.
   The existence of brake orbits in this case goes back to \cite{seifert_periodische_1948,bolotin_libration_1978}.
    Analogous to the conjecture of the existence of infinitely many periodic orbits for 
   Hamiltonian system on $(T^*M,\omega_0)$ (see e.g. \cite{long_multiple_2000,lu_conley_2009,mazzucchelli_lagrangian_2011}), 
   it is natural to consider the problem of the existence of infinitely many brake orbits
   for Hamilition system on $(T^*M,\omega_0)$. This problem has been studied in \cite{lu_existence_2007,lu_infinitely_2009} 
   for some special cases of $H$ and $M$. 
   The latest
   result in \cite{lu_conley_2009} gives the positive answer when $H$ is a $\tau$-periodic convex quadratic-growth Hamiltonian. More precisely,
   if there exist constants $0<h_1<h_2$, depending on the local coordiantes of $M$, such that 
\begin{equation}\label{sec1:01}
  \begin{aligned}
     &h_1|w|_q^2\leq\sum_{ij}\tfrac{\partial^2H}{\partial p_i\partial p_j}(t,q,p)w_iw_j\leq h_2|w|_q^2,\ \forall 0\neq w=\sum_{i=1}^N w_i\tfrac{\partial}{\partial p_i}\in T_q^*M,\\
&\left\lvert \tfrac{\partial^2H}{\partial q_i\partial p_j}(t,q,p)\right\rvert \leq h_2(1+|p|_q),\
 \left\lvert \tfrac{\partial^2H}{\partial q_i\partial q_j}(t,q,p)\right\rvert \leq h_2(1+|p|_q^2),
  \end{aligned}
\end{equation}
then the time-$\tau$ map of the Hamiltonian flow $\Psi_H^t$ (the Poincar$\acute{\text{e}}$ map) of $H$ has infinitely many periodic points sitting in 
   the zero section $M_0$ of $T^*M.$ Note that each $\tau$-periodic point of $\Psi_H^t$ sitting in $M_0$ corresponds to a 
   $\tau$-periodic brake orbit in $T^*M$.

    %We define a Hamiltonian function on a real symplectic manifold $(P,\omega,R)$ as a smooth function
    %$H: P\rightarrow\mathbb{R}$ satisfying $R^*H=H$. $H$ is called a time-dependent Hamiltonian if $H: \mathbb{R}
    %\times P\rightarrow\mathbb{R}$ is time-dependent.

    Let $\sigma$ be a closed 2-form (magnetic form) on the closed connected Riemannian manifold $M$. Equip $T^*M$ with the twisted symplectic
    form $\omega=\omega_0-\pi^*\sigma$, where $\pi: T^*M\rightarrow M$ is the natural projection. 
    %The anti-symplectic involution $R$ is given by $R(q,p)=(q,-p-2\theta(q))$ (see Theorem \ref{sec3:1}).
    When $\sigma=\text{d}\theta$ is exact, there exists an anti-symplectic involution 
    $R:T^*M\rightarrow T^*M$ defined by 
    \begin{equation*}
      R_1: (q,p)\mapsto (q,-p-2\theta(q))
    \end{equation*}
    (see in Section \ref{sec3:2}).
    We prove that 
    there are infinitely many brake orbits of the Hamiltonian system on $(T^*M,\omega)$ provided that 
    $H$ is a time-dependent Tonelli Hamilition and $R_1^*H(-t,\cdot)=H(t,\cdot)$. 
    We always assume that $H$ defines
     a global Hamiltonian flow on $T^*M.$  
     %Notice that a convex quadratic-growth Hamiltonian is a special case of a Tonelli Hamiltonian, 
    %we generize the result of \cite{luConleyConjectureHamiltonian2009}.

%   If we define a standard kinetic energy Hamiltonian $H$ on $T^*M,$ the Hamiltonian flow
%  of $H$ on $T^*M$, called a twisted geodesic flow, will be a physical model descrbing 
%the motion of a charge on $M$ in the magnetic field $\sigma.$ In this case, the existence of periodic orbits
%on fixed energy surface have been reserched by \dots. If $H$ is a Tonelli Hamiltonian, [Mazz]
%studied the multiplicity of periodic orbits on low energy surfaces.  

    %It is natural to consider brake orbits of a Hamiltonian 
    %system if we assume there exists an anti-symplectic
    %involution $R$ making $T^*M$ a real symplectic manifold. 
    %In this paper, we consider the Hamiltonian system of a Tonelli Hamiltonian $H$ on $(T^*M,\omega,R).$
    % We also assume the magnetic form $\sigma=\text{d}\theta$ on $M$ is exact, 
    %which means $\omega=\omega_0-\pi^*\text{d}\theta.$

    The main result is the following.
   
    \begin{thm}\label{sec1:th1}
      Let $M$ be a smooth closed manifold, $\sigma$ an exact magnetic form on $M$. 
      $H:S_1\times T^*M\rightarrow\mathbb{R}$ is a smooth 1-periodic Tonelli Hamilition on $(T^*M,\omega_0-\pi^*\sigma,R_1)$ 
    with global flow. 
    %and $a\in \mathbb{R}$ is a 
    %constant greater than
    % $$-\min_{q\in M}\left\{\int_{0}^1\min_{p\in T^*_qM}\{H(t,q,p)\}dt\right\}.$$
   %Assume that only finitely many 1-periodic brake orbits of the Hamilition system of $H$ have action less than $a$. 
   Then the Hamilition system of $H$ on $(T^*M,\omega_0-\pi^*\sigma)$ admits infinitely
   many brake orbits with period that is a power of $2$. 
    \end{thm}
    In particular, for the standard symplectic form $\omega_0$ on $T^*M$ and a Tonelli Hamiltonian $H$ on 
    $(T^*M,\omega_0,R_0)$, we have the following result which has been satated informally in \cite[Remark 6.1]{lu_corrigendum_2009}.
%\begin{cor} Let $(T^*M,\omega,R_1)$ and $H$ be defined as above. Let $\Phi_H^t$ be the global Hamiltonian flow of $H$
%  with respect to the twisted symplectic form $\omega.$ Then the Hamilition diffeomorphism $\Phi_H^1$ has 
%   infinitely many periodic points in ${Fix}(R_1)$.
%\end{cor}
\begin{cor}
  Let $M$ be a smooth closed manifold. 
      $H:S_1\times T^*M\rightarrow\mathbb{R}$ is a smooth 1-periodic Tonelli Hamilition on $(T^*M,\omega_0,R_0)$ 
    with global flow. 
    %and $a\in \mathbb{R}$ is a 
    %constant greater than
    % $$-\min_{q\in M}\left\{\int_{0}^1\min_{p\in T^*_qM}\{H(t,q,p)\}dt\right\}.$$
   %Assume that only finitely many 1-periodic brake orbits of the Hamilition system of $H$ have action less than $a$. 
   Then the Hamilition system of $H$ on $(T^*M,\omega_0)$ admits infinitely
   many brake orbits with period that is a power of $2$. 
\end{cor}

  Let $L: S_1\times TM\rightarrow\mathbb{R}$ be the dual Lagrangian of $H$. It is well known that the 
  Legendre transform sets up a one-to-one correspondence between the integer-periodic (brake) solutions of the Hamilition system
  of $H$ on $(T^*M,\omega_0-\pi^*\sigma)$ and integer-periodic (symmetric) solutions of the Euler-Lagrange system of $L+\theta$ on $M$, where $\theta$
  is a primitive of $\sigma$ (see in Section \ref{sec3:2}). 
  Therefore Theorem \ref{sec1:th1} can be equivalently stated in the Lagrangian formulation as follows.

  \begin{thm}\label{sec1:th3}
    Let $M$ be a smooth closed manifold with a smooth exact magnetic form $\sigma=\text{d}\theta$, 
    $L: S_1\times TM\rightarrow\mathbb{R}$ a smooth 1-periodic Tonelli Lagrangian with global flow and
    $(L+\theta)(-t,q,-v)=(L+\theta)(t,q,v)$. Then the Euler-Lagrangian system of $L+\theta$
    admits infinitely many symmetric periodic solutions with period that is a power of $2$.
  \end{thm}

   The proof of Theorem 1.3 is based on the method of \cite{long_multiple_2000}. More precisely, let $a$ be some 
   large number in $\mathbb{R}$, we assume by a 
   contradiction that there are only finite symmetric periodic solutions with mean Lagrangian action less than $a$, 
   then we deduce that there exists a symmetric periodic solution with the local Morse homology non-zero
   under infinitly many iterations. This will contradict the homological 
   vanishing result based on a time-reversible Bangert homotopy defined in Section \ref{sec7:5}.
   As far as we know, the Bangert homotopy defined in \cite{bangert_homology_1983} (see also
   \cite{lu_conley_2009}) may not be 
    a time-reversible one under any reparametrization. In Lemma \ref{sec7:3}, we give a construction of the time-reversible 
    Bangert homotopy, which is used to the proof of the homological vanishing theorem.

\section{Preliminaries}

\subsection{The symmetric loop space}\label{sec2:1}
We define $S_\tau:=\mathbb{R}/\tau\mathbb{Z},\, \forall \tau\in\mathbb{N}.$  Let $M$ be an $N$-dimensional 
smooth manifold with Riemannian metric $\langle\cdot,\cdot\rangle$. 
We denote by $H(\tau):=W^{1,2}(S_\tau,M)$ the space of absolutely continuous loops in $M$ with square-integrable weak 
derivative. 

Let $\gamma: S_\tau\rightarrow M$ be a $C^\infty$ $\tau$-periodic loop in $H(\tau)$, and
 $\gamma^*TM\rightarrow S_\tau$ be the pull-back of the  tangent bundle $TM$ by $\gamma.$
  We denote by $W^{1,2}(\gamma^*TM)$ the Hilbert space consisting of $W^{1,2}$-sections of the bundle 
$\gamma^*TM\rightarrow S_\tau$ with Hilbert structure given by 
\begin{equation}\label{sec2:11}
  \Angles{\xi}{\zeta}_\gamma:=\int_0^\tau [\angles{\xi(t)}{\zeta(t)}_{\gamma(t)}+\angles{\nabla_t\xi(t)}{\nabla_t\zeta(t)}_{\gamma(t)}]\text{d}t, \ \xi,\zeta\in W^{1,2}(\gamma^*TM),
\end{equation}
where $\nabla_t$ denotes the covariant derivative with respect to the Levi-Cicita connection on the Riemannian manifold $(M,\angles{\cdot}{\cdot})$.
 It is well known that $H(\tau)$ is a Hilbert manifold and the tangent space 
at $\gamma$ is $W^{1,2}(\gamma^*TM)$, and \eqref{sec2:11} defines a Hilbert-Riemannian metric on $H(\tau)$.
In this case, $H(\tau)$ turns out to be a complete Hilbert-Riemannian manifold (see \cite{klingenberg_lectures_1978}).

If $\gamma\in H(\tau)$ is time-reversible, i.e., $\gamma(-t)=\gamma(t), \forall t\in\mathbb{R},$
we call $\gamma$ a symmetric loop. We follow the notation in \cite{lu_conley_2009} and denote by $EH(\tau)$ the symmetric (even) 
loop space consisting of time-reversible loops in $H(\tau)$.

%There is a natural manifold structure of $EH(\tau)$ such that $EH(\tau)$ is a submanifold of $H(\tau)$. 
Now we assume $\gamma\in EH(\tau)$ is smooth. Let $EW^{1,2}(\gamma^*TM)$ be the subspace of $W^{1,2}(\gamma^*TM)$
consisting of time-reversible $W^{1,2}$-sections. Then $EW^{1,2}(\gamma^*TM)$ is a
 Hilbert subspace of $ W^{1,2}(\gamma^*TM)$ with the Hilbert structure defined in \eqref{sec2:11}.
Let $\epsilon>0$ be a constant smaller than
the injectivity radius of $(M,\angles{\cdot}{\cdot})$. We define $N_\epsilon:=\{v\in T_qM\,|\, q\in M, |v|_q<\epsilon\}$
and $\gamma^*N_\epsilon:=\gamma^*TM\cap N_\epsilon$. Then $EW^{1,2}(\gamma^*N_\epsilon)\subset EW^{1,2}(\gamma^*TM)$ is an open set of sections along $\gamma$ that 
take values inside $N_\epsilon$.   
we define a injective map
\begin{equation*}
  \exp_\gamma: EW^{1,2}(\gamma^*N_\epsilon)\rightarrow H(\tau)
\end{equation*}
as 
$$(\exp_\gamma\xi)(t):=\exp_{\gamma(t)}(\xi(t)),\quad \forall\xi\in EW^{1,2}(\gamma^*N_\epsilon).$$
It is easy to know that $(\exp_\gamma\xi)(-t)=(\exp_\gamma\xi)(t).$ Let $\mathcal{U}_\gamma$ be the image of $\exp_\gamma$,
 the following natural atlas
 $$\left\{(\exp_\gamma^{-1},\mathcal{U}_\gamma),\quad\text{for time-reversible}\ \gamma\in C^\infty(S_\tau,M)\right\}$$ 
define the differentiable structure of $EH(\tau)$. If we restrict $\Angles{\cdot}{\cdot}$ 
 on $EH(\tau)$, \eqref{sec2:11} also defines a Hilbert-Riemannian metric on $EH(\tau)$
and makes the latter into a complete Hilbert-Riemannian manifold.

 For each $\gamma\in EH(\tau)$, we define a homotopy 
$\gamma_s$ connecting $\gamma_0:=\gamma(0)$ and $\gamma_1:=\gamma$ as
\begin{equation*}
  \gamma_s(t):=
   \begin{cases}
    \gamma(st), &t\in [0,\tfrac{\tau}{2})\\
    \gamma(s(\tau-t)), & t\in[\tfrac{\tau}{2},\tau].
  \end{cases}
\end{equation*} 
Therefore $\gamma$ is contractible and $EH(\tau)$ is a subspace of the contractible loop space of $H(\tau).$

Let $S_{n\tau}:=\mathbb{R}/n\tau\mathbb{Z}$ be the $n$-fold covering of $S_\tau.$ We define the iteration map 
\begin{equation*}
  \psi^{[n]}: H(\tau)\hookrightarrow H(n\tau)
\end{equation*}
such that $\gamma^{[n]}:=\psi^{[n]}(\gamma)$ is defined as the composition of $\gamma$ with the $n$-fold covering map of $S_\tau.$
For each $\gamma\in H(\tau)$, it is easy to know that the differential $\text{d}\psi^{[n]}$ of $\psi^{[n]}$ defines an iteration map 
\begin{equation*}
  \text{d}\psi^{[n]}: W^{1,2}(\gamma^*TM)\hookrightarrow W^{1,2}(\gamma^{[n]*}TM),
\end{equation*}
that is similar to that of $\psi^{[n]}$, and $\exp_{\gamma^{[n]}}\circ\text{d}\psi^{[n]}=\psi^{[n]}\circ \exp_\gamma.$ 
We notice that both the maps $\psi^{[n]}$ and $\text{d}\psi^{[n]}$ restrict as iteration maps 
on $EH(\tau)$ and $EW^{1,2}(\gamma^*TM)$ respectively.

\subsection{Hamiltonian settings and Legendre transform}\label{sec2:2}

Let $H$ be a smooth $\tau$-periodic Tonelli Hamiltonian on symplectic manifold
$(T^*M,\omega_0,R_0)$.
The Tonelli conditions (H1)(H2) admits the global existence of the inverse Legendre teansform, given by
\begin{equation*}
  \mathcal{L}_H: S_\tau\times T^*M\rightarrow S_\tau\times TM,\; (t,q,p)\rightarrow(t,q,\partial_pH(t,q,p)).
\end{equation*}
$\mathcal{L}_H$ defines a fiber-preserving $C^\infty$-diffeomorphism. 

In this case, the Fenchel transform of the Tonelli Hamiltonian $H: S_1\times T^*M\rightarrow\mathbb{R}$, defined by 
\begin{equation*}
   L(t,q,v):=\max\{p[ v]-H(t,q,p)\, |\, p\in T_q^*M\},\quad  \forall (t,q,v)\in S_1\times TM,
\end{equation*}
 will be a $C^\infty$ Lagrangian function and satisfies
\begin{equation*}
  L(t,q,v)=p(t,q,v)[v]-H(t,q,p(t,q,v)),
\end{equation*}
where $p=p(t,q,v)$ is the unique point determined by the equality $v=\partial_pH(t,q,p)$. We refer the reader to 
\cite{fathi_weak_2008} for more details about Tonelli functions and the 
properties of  
 Legendre transform.
It is easy to know that the dual Lagrangian $L$  satisfies the following properties.

(L1) the fiberwise Hessian of $L$ is positive definite, i.e. 
\begin{equation*}
  \sum_{i,j=1}^N\frac{\partial^2L}{\partial v_i\partial v_j}(t,q,v)u_iu_j>0
\end{equation*}
 for all $(t,q,v)\in S_\tau\times TM$ and nonzero $ u=\sum_{i=1}^N u_i\tfrac{\partial}{\partial v_i}\in T_qM.$

 (L2) $L$ is fiberwise superlinear, i.e.
\begin{equation*}
  \lim_{|v|_q\rightarrow\infty}\frac{L(t,q,v)}{|v|_q}=\infty,
\end{equation*}
for all $(t,q)\in S_\tau\times M.$

(L3) $L(-t,q,-v)=L(t,q,v),\, \forall(t,q,v)\in S_\tau\times TM.$

(L1)(L2) are the Tonelli conditions of $L$ and we say $L$ the Tonelli Lagrangian.
We define the Euler-Lagrange system of $L$ on the configuration space $M$ by 
\begin{equation}\label{sec2:22}
  \frac{\text{d}}{\text{d}t}\frac{\partial L}{\partial v_j}(t,\gamma(t),\dot{\gamma}(t))=\frac{\partial L}{\partial q_j}(t,\gamma(t),\dot{\gamma}(t)),\; j=1,\dots,N.
\end{equation} 
in any local coordinates $(q_1,\dots,q_N).$ It is well known that
the Euler-Lagrange flow $\Psi_L^t$ is conjugated by $\mathcal{L}_H$ to the Hamiltonian flow $\Psi_H^t$ of $H$ (see \cite[Theorem 3.4.2]{fathi_weak_2008}).
Therefore, there is a one-to-one correspondence between the periodic solutions of the Euler-Lagrange system of $L$
and the periodic solutions of the Hamiltonian system of $H$. In particular, $\gamma(t)$ is a $\tau$-periodic symmetric solution of \eqref{sec2:22}  
if and only if $x(t)=(\gamma(t),\partial_vL(t,\gamma(t),\dot{\gamma}(t)))$ is a $\tau$-periodic brake orbit of
the Hamilition system \eqref{sec1:02}
with $H$ satisfying (H1)(H2).
%It is well known that a curve $\mathbb{R}\rightarrow T^*M, t\mapsto x(t):=(\gamma(t),\gamma^*(t))$ is a solution of the Hamiltonian 
%system of $H$ if and only if 
%$$\gamma^*(t)=\partial_vL(t,\gamma(t),\dot{\gamma}(t)),\quad \forall t\in\mathbb{R}$$
% and $\gamma$ is a solution
%of the Euler-Lagrangian system on $M$,

%Besides, $x(t)=(\gamma(t),\gamma^*(t))$ is a brake orbit of $H$ if and only if 
%$\gamma(-t)=\gamma(t)$ is time-reversible and $\gamma(t+\tau)=\gamma(t).$ 

For each $\tau\in \mathbb{N}$, we define the mean action functional with respect to $L$ on $H(\tau)$ by 
\begin{equation*}
  A^{[\tau]}(\gamma)=\frac{1}{\tau}\int_0^\tau L(t,\gamma(t),\dot{\gamma}(t))\text{d}t.
\end{equation*}
We denote by $EA^{[\tau]}:=A^{[\tau]}|_{EH(\tau)}.$ 
$A^{[\tau]}$ and $EA^{[\tau]}$ are Gateaux-differentiable on $H(\tau)$ and $EH(\tau)$ respectively.
Each $\tau$-periodic symmetric solution $\gamma:S_\tau\rightarrow M$ of the Euler-Lagrange system
\eqref{sec2:22} is smooth and an extremal of $A^{[\tau]}$ and $EA^{[\tau]}$ (See e.g.\cite[Corollary 2.2.11]{fathi_weak_2008}).

Let $\gamma\in EH(\tau)$ be the weak $C^0$-extremal of $EA^{[\tau]}$, 
i.e. $\gamma$ satisfies $dEA^{[\tau]}(\gamma)[\xi]=0$ for each smooth section 
$\xi\in EW^{1,2}(\gamma^*TM)$. Denote by 
$$\Xi: W^{1,2}(\gamma^*TM)\rightarrow W^{1,2}(\gamma^*TM),\  \xi(t)\mapsto\xi(-t)$$ the involution,
we have 
%$A^{[\tau]}(\Xi(\xi))=A^{[\tau]}(\xi)$ and 
\begin{equation}\label{sec2:23}
  dA^{[\tau]}(\gamma)[\Xi(\xi)]=dA^{[\tau]}(\gamma)[\xi].
\end{equation}
We rewrite $\xi=\tfrac{1}{2}(\xi+\Xi(\xi))+\tfrac{1}{2}(\xi-\Xi(\xi)):=\xi_0+\xi_1.$ Since $\xi_0\in EW^{1,2}(\gamma^*TM)$ 
is time-reversible, we have $dA^{[\tau]}(\gamma)[\xi_0]=dEA^{[\tau]}(\gamma)[\xi_0]=0$ and 
\begin{equation}\label{sec2:24}
  \begin{aligned}
  dA^{[\tau]}(\gamma)[\xi]&= dA^{[\tau]}(\gamma)[\xi_1]\\
                          &=\tfrac{1}{2}dA^{[\tau]}(\gamma)[\xi]-\tfrac{1}{2}dA^{[\tau]}(\gamma)[\Xi(\xi)]\\
                          &=0\quad \text{by}\ \eqref{sec2:23}
\end{aligned} 
\end{equation}
for each smooth section $\xi\in W^{1,2}(\gamma^*TM).$ Therefore, $\gamma$ is a weak $C^0$-extremal of $A^{[\tau]}$. However, $\gamma$ may not be a solution of the 
Euler-Lagrange equation \eqref{sec2:22} unless $\gamma$ is $C^1$.
 
It is noted that each Tonelli Lagrangian can be modified to be quadratically growth in $v$ for $|v|$ large
(see in Section 5 or \cite{abbondandolo_high_2007}). We assume there exist $0<l_1<l_2$ depending on the
local coordinates of $M$ such that
\begin{equation}\label{sec2:311}
  \begin{aligned}
  &\sum_{ij}\tfrac{\partial^2L}{\partial v_i\partial v_j}(t,q,v)u_iu_j\geq l_1|u|_q^2,\ \forall\ 0\neq u=\sum_{i=1}^N u_i\tfrac{\partial}{\partial v_i}\in T_qM,\\
&\left\lvert \tfrac{\partial^2L}{\partial v\partial v}(t,q,v) \right\lvert\leq l_2,\
\left\lvert \tfrac{\partial^2L}{\partial q\partial v}(t,q,v)\right\rvert \leq l_2(1+|v|_q),\
\left\lvert \tfrac{\partial^2L}{\partial q\partial q}(t,q,v)\right\rvert \leq l_2(1+|v|_q^2).
\end{aligned} 
\end{equation}
 In this case, $L$ is the so-called convex quadratic-growth Lagrangian.
 We introduce the properties of the mean action functional 
 $A^{[\tau]}: H(\tau)\rightarrow \mathbb{R}$ of $L$ as follows.

(i) $A^{[\tau]}$ is $C^1$-continuous and twice Gateaux-differentiable;

(ii) Critical points are precisely the (smooth) solutions of the Euler-Lagrange system of $L$;

 (iii) $A^{[\tau]}$ satisfies the Palais-Smale condition;

 See \cite[Proposition 3.1]{abbondandolo_smooth_2009-3} and \cite[Proposition 4.2]{abbondandolo_high_2007} for details. 
 There are similar properties for $EA^{[\tau]}$. More precisely,

(i') since $EA^{[\tau]}$ is the restriction of $A^{[\tau]}$ on $EH(\tau)$, $EA^{[\tau]}$ 
is $C^1$-continuous and twice Gateaux-differentiable.

(ii') by \eqref{sec2:24}, each critical point $\gamma$ of $EA^{[\tau]}$ is a critical point of $A^{[\tau]}$.
Therefore $\gamma$ is smooth and is a symmetric solution of the Euler-Lagrange system of $L$.

(iii') let $\{\gamma_n\}\subset EH(\tau)$ be a sequence such that $EA^{[\tau]}(\gamma_n)$ is bounded and 
\begin{equation*}
  \| dEA(\gamma_n)|_{EW^{1,2}(\gamma_n^*TM)}\|_*\rightarrow 0,\quad\text{ as }n\rightarrow \infty.
\end{equation*}
%Here $T_{\gamma_n}EW^{1,2}(S_\tau,M)$ denotes the tangent space at some smooth loop $\gamma_n'$,
%which has a coordinate neighborhood containing $\gamma_n$, and $\|\cdot \|_*$ denotes the dual norm on $T^*_{\gamma_n}EW^{1,2}(S_\tau,M)$.
Here $\|\cdot\|_*$ is the dual norm of the dual space $EW^{1,2}(\gamma_n^*TM)^*$. We have
 \begin{align*}
   &\| dA(\gamma_n)|_{W^{1,2}(\gamma_n^*TM)}\|_*\\
  =&\sup_{\|\xi\|=1}|dA^{[\tau]}(\gamma_n)[\xi]|
  =\sup_{\|\xi\|=1}|\tfrac{1}{2}dEA^{[\tau]}(\gamma_n)[\xi_0]|\\
  =&\sup_{\|\xi_0\|=1}|dEA^{[\tau]}(\gamma_n)[\xi_0]|\rightarrow 0 \text{ as }n\rightarrow \infty.
\end{align*}
Therefore $\{\gamma_n\}$ is a Palais-Smale sequence of $A^{[\tau]}$, which is compact in $H(\tau).$
 We rewrite $\{\gamma_n\}$ as the 
subsequence of $\{\gamma_n\}$ such that $\gamma_n\rightarrow \gamma$ in $H(\tau).$
$\gamma$ is the uniform limit of $\{\gamma_n\},$ thus $\gamma(\tau-t)=\gamma(t).$
We conclude that $EA^{[\tau]}$ satisfies the Palais-Smale condition.

 We will simply write $A=A^{[1]}$ and $EA=EA^{[1]}.$
 Let $\psi^{[n]}$ be an iteration map defined respectively on $H(\tau)$ and $EH(\tau),$ thus
$$A^{[n\tau]}\circ\psi^{[n]}=A^{[\tau]},\quad EA^{[n\tau]}\circ\psi^{[n]}=EA^{[\tau]}.$$

\subsection{Morse index and Maslov-type index }\label{sec2:3}

Let $L: S_\tau\times TM\rightarrow \mathbb{R}$ be a smooth time-reversible convex quadratic-growth Lagrangian. 
 We introduce the Morse index and nullity of $A^{[\tau]}$ and $EA^{[\tau]}$ at $\gamma$.
 
Since $\gamma$ is contractible, $\gamma^*TM$ is a trivial bundle, and we can choose a smooth orthogonal trivialization
\begin{equation*}
  S_\tau\times \mathbb{R}^N\rightarrow \gamma^*TM,\quad (t,q)\mapsto\Phi(t)q
\end{equation*}
with $\Phi(-t)=\Phi(t),\ \forall t\in S_\tau.$

Let $\rho$ be a positive number smaller than the injectivity
radius of $(M,\angles{\cdot}{\cdot})$ and $B_\rho:=B_\rho^N(0)$ be an open ball in $\mathbb{R}^N$ centered at 0 with radius $\rho$. 
Then for each $k\in\mathbb{N}$, the following differentiable injection gives a coordinate chart on $H(k\tau)$ containing $\gamma^{[k]}$,
\begin{equation}\label{sec2:34}
  \Theta^{[k]}_\gamma: W^{1,2}(S_{k\tau},B_\rho)\rightarrow \mathcal{U}_{\gamma^{[k]}}\subset H(k\tau),\ \xi\mapsto \exp_{\gamma^{[k]}}(\Phi^{[k]}\xi).
\end{equation}
Here $\Phi^{[k]}$ is the composition of $\Phi$ with the $k$-fold covering map of $S_\tau.$ 
Therefore, in a neighborhood of $\gamma^{[k]}$ inside $\mathcal{U}_{\gamma^{[k]}}$, $A^{[k\tau]}$ restricts 
as a functional $ W^{1,2}(S_{k\tau},B_\rho)\rightarrow\mathbb{R},$ and
\begin{align*}
&A^{[k\tau]}\circ\Theta_\gamma^{[k]}(\xi)\\
=&\frac{1}{k\tau}\int_0^{k\tau} L\left(t,\Theta_\gamma^{[k]}(\xi)(t),\frac{\text{d}}{\text{d}t}\Theta^{[k]}(\gamma^{[k]})(t)\right)\text{d}t\\
                                         =&\frac{1}{k\tau}\int_0^{k\tau} L\left(t,\Theta_\gamma^{[k]}(\xi)(t),\text{d}(\Theta_\gamma^{[k]})(\xi)\cdot\dot{\xi}+\frac{\text{d}}{\text{d}t}\Theta_\gamma^{[k]}(t)\circ\xi\right)\text{d}t.
\end{align*}

Now we define $\pi^{[k]}_{\gamma}: S_{\tau}\times B_{\rho}\times\mathbb{R}^N\rightarrow S_\tau\times TM$ as
\begin{equation*}
  (t,q,v)\mapsto\left(t,\Theta_\gamma^{[k]}(q)(t),\text{d}(\Theta_\gamma^{[k]})(q)\cdot v+\frac{\text{d}}{\text{d}t}\Theta_\gamma^{[k]}(t)\circ q\right)
\end{equation*}
then 
\begin{equation*}
  A^{[k\tau]}\circ\Theta_\gamma^{[k]}(\xi)=\frac{1}{k\tau}\int_0^{k\tau}L\circ\pi^{[k]}_{\gamma}(t,\xi(t),\dot{\xi}(t))\text{d}t,
\end{equation*}
and $\mathbf{0}=(\Theta_\gamma^{[k]})^{-1}(\gamma^{[k]})\in W^{1,2}(S_{k\tau},B_{\rho})$ is a critical point of $A^{[k\tau]}\circ\Theta_\gamma^{[k]}.$ 

Since $\gamma\in EH(\tau)$ is symmetric, and $\Phi(-t)=\Phi(t)$, we have 
\begin{equation*}
  \Theta^{[k]}_\gamma(\xi)(-t)=\Theta^{[k]}_\gamma(\xi)(t),\;\;\forall \xi\in EW^{1,2}(S_{k\tau},B_\rho).
\end{equation*} 
Thus
 $\Theta^{[k]}_\gamma$ restricts as a coordinate chart of $EH(k\tau)$ at $\gamma.$ 
Besides, we compute that  $\forall (t,q,v)\in S_\tau\times B_\rho\times\mathbb{R}^N,$
\begin{align*}
  L\circ\pi^{[k]}_{\gamma}(-t,q,-v)=L\circ\pi^{[k]}_{\gamma}(t,q,v),\\
   L\circ\pi^{[k]}_{\gamma}(t+1,q,v)=L\circ\pi^{[k]}_{\gamma}(t,q,v).
\end{align*}
 Therefore $L\circ\pi^{[k]}_{\gamma}$ satisfies the property \eqref{sec2:311} 
 and is a time-reversible convex quadratic-growth Lagrangian.
%We have
%\begin{equation*}
%  EA^{[k\tau]}\circ\Theta_\gamma^{[k]}(\xi)=\frac{1}{k\tau}\int_0^{k\tau}L\circ\pi^{[k]}_{\gamma}(t,\xi(t),\dot{\xi}(t))\text{d}t,
%\end{equation*}
%and $\mathbf{0}=(\Theta_\gamma^{[k]})^{-1}(\gamma^{[k]})\in EW^{1,2}(S_{k\tau},B_{\rho})$ is a critical point of $EA^{[k\tau]}\circ\Theta_\gamma^{[k]}.$ 

In the following, we will simply write $L$, $A^{[k\tau]}$, $EA^{[k\tau]}$ for
$L\circ\pi^{[k]}_{\gamma}$, $A^{[k\tau]}\circ\Theta_\gamma^{[k]}$ and $EA^{[k\tau]}\circ\Theta_\gamma^{[k]}$ respectively.
Then $\gamma^{[k]}$ will be identified with the point $\mathbf{0}$ in the Hilbert space $W^{1,2}(S_{k\tau},\mathbb{R}^N).$ 
%and as the critical point of $A^{[\tau]}$ and $EA^{[\tau]}$ respectively on $H(\tau)$ and $EH(\tau).$
Without loss of generality, we assume $M=\mathbb{R}^N$ and
\begin{align*}
  & L: S_{k\tau}\times \mathbb{R}^{2N}\rightarrow\mathbb{R},\\
  &A^{[k\tau]}: W^{1,2}(S_{k\tau},\mathbb{R}^N)\rightarrow\mathbb{R},\\
  &EA^{[k\tau]}: EW^{1,2}(S_{k\tau},\mathbb{R}^N)\rightarrow\mathbb{R}.
\end{align*}

We denote by {Hess$A^{[k\tau]}(\mathbf{0})$ the Hessian of $A^{[k\tau]}$ at the critical point $\mathbf{0}\in W^{1,2}(S_{k\tau},\mathbb{R}^N)$.}
For each $ \xi,\eta\in W^{1,2}(S_{k\tau},\mathbb{R}^N)$, we have
\begin{equation*}
  \text{Hess}A^{[k\tau]}(\mathbf{0})[\xi,\eta]=\frac{1}{k\tau}\int_0^{k\tau} (P(t)\dot{\xi}\cdot \dot{\eta}+Q(t)\xi\cdot\dot{\eta}+Q^T(t)\dot{\xi}\cdot\eta+R(t)\xi\cdot\eta)\text{d}t,
\end{equation*} 
 where 
\begin{equation}\label{sec2:35}
  \begin{aligned}
     P_{i,j}(t)&:=\frac{\partial^2L}{\partial v_i\partial v_j}(t,0,0),\\
      Q_{i,j}(t)&:=\frac{\partial^2L}{\partial q_i\partial v_j}(t,0,0),\\
       R_{i,j}(t)&:=\frac{\partial^2L}{\partial q_i\partial q_j}(t,0,0).
  \end{aligned}
 \end{equation}
 and $P(1-t)=P(t),\,Q(1-t)=-Q(t),\,R(1-t)=R(t).$ $Q^T(t)$ is the transpose of $Q(t)$ for each $t\in S_{k\tau}.$
 Here we do not distinguish $P,Q,R$ with their $k$-iterations. By the definition of $EA^{[k\tau]},$ we have
 \begin{equation*}
   \text{Hess}EA^{[k\tau]}(\mathbf{0})=\text{Hess}A^{[k\tau]}(\mathbf{0})|_{EW^{1,2}\times EW^{1,2}}.
 \end{equation*} 
 
 Let 
 \begin{align*}
   W^{1,2}(S_{k\tau},\mathbb{R}^N)&=W^{1,2}(S_{k\tau},\mathbb{R}^N)^+\oplus W^{1,2}(S_{k\tau},\mathbb{R}^N)^0\oplus W^{1,2}(S_{k\tau},\mathbb{R}^N)^-,\\
   EW^{1,2}(S_{k\tau},\mathbb{R}^N)&=EW^{1,2}(S_{k\tau},\mathbb{R}^N)^+\oplus EW^{1,2}(S_{k\tau},\mathbb{R}^N)^0\oplus EW^{1,2}(S_{k\tau},\mathbb{R}^N)^-
 \end{align*}
 be respectively the Hess$A^{[k\tau]}(\mathbf{0})$-orthogonal decompositions according to $\text{Hess}A^{[k\tau]}(\mathbf{0})$ and 
 $\text{Hess}EA^{[k\tau]}(\mathbf{0})$ being positive, null, and negative definite. 
Then the Morse index and nullity of $\gamma$ for $A^{[k\tau]}$ are defined by
\begin{align*}
  m^-(A^{[k\tau]},\gamma)&:=\dim W^{1,2}(S_{k\tau},\mathbb{R}^N)^-,\\
  m^0(A^{[k\tau]},\gamma)&:=\dim W^{1,2}(S_{k\tau},\mathbb{R}^N)^0,
\end{align*}
 which are well-defined (may be infinity) since they do not depend on the choice of the trivialization of $\gamma^{[k]*}TM$ and the coordinate chart $\Theta_\gamma^{[k]}.$
 Similarly, 
 \begin{align*}
  m^-(EA^{[k\tau]},\gamma)&:=\dim EW^{1,2}(S_{k\tau},\mathbb{R}^N)^-,\\
  m^0(EA^{[k\tau]},\gamma)&:=\dim EW^{1,2}(S_{k\tau},\mathbb{R}^N)^0
\end{align*}
 are defined as the Morse index and nullity of $\gamma^{[k]}$ for $EA^{[k\tau]}$. 
 It is well known that for each nonzero $k\in\mathbb{N}$, every critical point of $A^{[k\tau]}$ and 
 $EA^{[k\tau]}$ has finite Morse index and nullity (see \cite{benci_periodic_1986}), and
 \begin{align}
   &0\leq m^-(EA^{[k\tau]},\gamma)\leq m^-(A^{[k\tau]},\gamma),\label{sec2:32}\\
   &0\leq m^0(EA^{[k\tau]},\gamma)\leq m^0(A^{[k\tau]},\gamma)\leq 2N\label{sec2:31}
 \end{align}
 (see \cite[Theorem 3.4]{lu_infinitely_2009}). 

 On the other hand, there is a Maslov-type index related to $\gamma.$ 
 Recall the definition of Legendre transform in Section \ref{sec2:2}, 
 we have the dual Hamiltonian $H: S_\tau\times \mathbb{R}^{2N}\rightarrow\mathbb{R}$ defined by
 \begin{align*}
   H(t,q,p)&=p[v(t,q,p)]-L(t,q,v(t,q,p)),\\
   p&=\partial_vL(t,q,v(t,q,p)).
 \end{align*}
 The corresponding Hamiltonian system is 
 \begin{align*}
   \dot{p}&=\partial_qH(t,q,p)=-\partial_qL(t,q,v(t,q,p)),\\
   \dot{q}&=\partial_pH(t,q,p)=v(t,q,p).
 \end{align*}
and loop $x(t)=(0,\partial_vL(t,0,0))$ is the corresponding $\tau$-periodic solution.
Linearization along $x$ yields the linear periodic Hamiltonian system
\begin{equation}\label{sec2:310}
  \begin{aligned}
  \dot{\xi}&=Q^TP^{-1}\xi+Ry-Q^TP^{-1}Qy,\\
  \dot{y}&=P^{-1}\xi-P^{-1}Qy.
\end{aligned}
\end{equation}

We denote $\mathbf{u}=\left(\begin{matrix}
  \xi\\
  y
\end{matrix}\right)$ and  
\begin{equation*}
  B(t)=\left(\begin{matrix}
    P(t)^{-1}&-P(t)^{-1}Q(t)\\
    -Q(t)^TP(t)^{-1}&Q(t)^TP(t)^{-1}Q(t)-R(t)
  \end{matrix}\right)
\end{equation*}
Then \eqref{sec2:310} can be simply written as 
\begin{equation*}
  \dot{\mathbf{u}}(t)=JB(t)\mathbf{u},
\end{equation*}
Here $J=\left(\begin{matrix}
  0&-I\\
  I&0
\end{matrix}\right).$
The fundamental solution $\Psi_\gamma$ of this linear Hamiltonian system with initional point $\Psi_\gamma(0)=I$
is a symplectic path along $ x(t)$. Here the subscript $\gamma$ implies that the symplectic path depends on $\gamma$,
since the Hamilition loop $x(t)$ is essentially the Legendre transform of the Lagrangian loop $\gamma(t)$. 
The Maslov-type index pair and the $k$-iteration index pair of $\Psi_\gamma$
\begin{equation*}
  (i(\Psi_\gamma),\nu(\Psi_\gamma)),\quad (i(\Psi_\gamma,k),\nu(\Psi_\gamma,k))
\end{equation*} are defined by C. Conley, E. Zehnder, Y. Long and C. Viterbo (see e.g. \cite{long_index_2002}).
 It is well known that the Morse index pair $(m^{-}(A^{[k\tau]},\gamma),m^0(A^{[k\tau]},\gamma))$
of the functional $A^{[k\tau]}$ at $\gamma^{[k]}$ coincides with the $k$-iteration Maslov-type index of $\Psi_\gamma$, i.e
\begin{equation*}
  (m^-(A^{[k\tau]},\gamma),m^0(A^{[k\tau]},\gamma))=(i(\Psi_\gamma,k),\nu(\Psi_\gamma,k))
\end{equation*}
(See \cite{viterbo_new_1990,long_indexing_1998,abbondandolo_morse_2003}).

Besides, $x(-t)=R_0x(t)$ implies that $x:S_\tau\rightarrow\mathbb{R}^{2N}$ is in fact a brake orbit. Then the 
Maslov-type $L_0$-index pair and the corresponding $k$-iteration index pair of $\Psi_\gamma$
\begin{equation*}
  (i_{L_0}(\Psi_\gamma),\nu_{L_0}(\Psi_\gamma)),\quad (i_{L_0}(\Psi_\gamma,k),\nu_{L_0}(\Psi_\gamma,k))
\end{equation*}
defined in \cite{liu_maslov-type_2007} applies. By \cite[Theorem 3.4]{lu_infinitely_2009} and \cite[Proposition 6.1]{liu_iteration_2014},
 we have the similar result for the relation between the Morse index pair $(m^{-}(EA^{[\tau]},\gamma),m^0(EA^{[\tau]},\gamma))$ 
and the Maslov-type $L_0$-index pair, that is 
\begin{equation*}
  (m^-(EA^{[k\tau]},\gamma),m^0(EA^{[k\tau]},\gamma))= (i_{L_0}(\Psi_\gamma,k)+N,\nu_{L_0}(\Psi_\gamma,k)).
\end{equation*}

Let 
\begin{equation*}
  \widehat{m}^-(A^{[\tau]},\gamma):=\lim_{k\rightarrow\infty}\frac{m^-(A^{[k\tau]},\gamma)}{k},\quad 
  \widehat{m}^-(EA^{[\tau]},\gamma):=\lim_{k\rightarrow\infty}\frac{m^-(EA^{[k\tau]},\gamma)}{k}
\end{equation*}
 be the mean Morse index of $m^-(A^{[k\tau]},\gamma)$  and 
$m^-(EA^{[k\tau]},\gamma)$ respectively. We have 
\begin{equation*}
  \widehat{m}^-(EA^{[\tau]},\gamma)=\widehat{i}_{L_0}(\Psi_\gamma)=\frac{1}{2}\widehat{i}(\Psi_\gamma)=\frac{1}{2}\widehat{m}^-(A^{[\tau]},\gamma),
\end{equation*}
where 
\begin{equation*}
  \widehat{i}_{L_0}(\Psi_\gamma):=\lim_{k\rightarrow\infty}\frac{i_{L_0}(\Psi_\gamma,k)}{k},\quad
  \widehat{i}(\Psi_\gamma):=\lim_{k\rightarrow\infty}\frac{i(\Psi_\gamma,k)}{k},
\end{equation*}
 and the middle equality is quoted from \cite[Proposition 4.1]{liu_iteration_2014}.
Since 
\begin{equation*}
  i(\Psi_\gamma,k)+\nu(\Psi_\gamma,k)\leq k\widehat{i}(\Psi_\gamma)+N
\end{equation*}
(see e.g. \cite[p. 213 Theorem 2]{long_index_2002}), we have 
\begin{equation}\label{sec2:33}
  0\leq m^-(EA^{[k\tau]},\gamma)+m^0(EA^{[k\tau]},\gamma)\leq N,\text{ if }  \widehat{m}^-(EA^{[\tau]},\gamma)=0.
\end{equation}

\section{The Hamilition flow on \texorpdfstring{$(T^*M,\omega)$}{} and the dual Lagrangian flow }\label{sec3:2}
As before, we assume $M$ is a closed connected smooth manifold. 
Let $\sigma$ be a smooth magnetic form on $M$, and $\omega=\omega_0-\pi^*\sigma$ be a twisted symplectic form on $T^*M$.
%\subsection{"Untwist" the symplectic manifold \texorpdfstring{$(T^*M,\omega)$}{}}
If $\sigma$ is exact, we show that there exists a symplectomorphism between the symplectic manifolds
 $(T^*M,\omega)$ and $(T^*M,\omega_0)$.
%  Let $R: T^*M\rightarrow T^*M $ be given by $ R(q,p)= (q,-p-2\theta(q)),\, \forall (q,p)\in T^*M$. 
%  Thus $R^2=id$ and $R^*\omega=-\omega$. We give a proof for the existence of a symplectomorphism of real symplectic manifolds $(T^*M,\omega,R)$
%  and $(T^*M,\omega_0,R_0)$ which commutes with the corresponding anti-symplectic involutions. 

\begin{thm}\label{sec3:1}
   %Let $(T^*M,\omega,R)$ be a real symplectic manifold, where $\omega=\omega_0-\pi^*\text{d}\theta$. 
   If $\sigma=\text{d}\theta$ is exact, there exists a global diffeomorphism $\Phi:(T^*M,\omega_0)\rightarrow (T^*M,\omega)$
such that $\Phi$ is symplectic, i.e. $\Phi^*\omega=\omega_0$. 
%Therefore, $R:=\Phi\circ R_0\circ\Phi^{-1}$ is a natural anti-symplectic involution on $(T^*M,\omega).$
\end{thm}
%(2) the following diagram commutes. 
%\begin{center}
%  \begin{tikzcd}
%(T^*M,\omega_0) \arrow[r,"R_0"] \arrow[d,"\Phi"] & (T^*M,\omega_0) \arrow[d,"\Phi"]\\
%(T^*M,\omega) \arrow[r,"R"]             & (T^*M,\omega) 
%\end{tikzcd}
%\end{center}
%

~\\
\Proof.  We employ the so-called deformation method of J. Moser. Define 
$$\omega_t=\omega_0-t\pi^*\text{d}\theta, \quad 0\leq t\leq1,$$
such that $\omega_1=\omega.$ Note that $\omega_t$ is 
a symplectic form for each $t\in[0,1].$ We look for a whole family $\phi^t: T^*M\rightarrow T^*M$ of diffeomorphism satisfying
$ \phi^0=\text{id}$ and
 $$(\phi^t)^*\omega_t=\omega_0, \quad 0\leq t\leq1.$$
We denote by $\Phi=\phi^1$. Then $\Phi$ will be the solution to our problem.  
 
 We shall construct a $t$-dependent vector field $X_t$ generating $\phi^t$ as its flow.
 Since
 \begin{align*}
   0=\frac{\text{d}}{\text{d}t}(\phi^t)^*\omega_t=(\phi^t)^*\text{d}\left(i_{X_t}\omega_t-\pi^*\theta\right).
 \end{align*}
 $X_t$ has to satisfy the linear equation
$$i_{X_t}\omega_t=\pi^*\theta.$$
Note that $\omega_t$ are nondegenerate for $0\leq t\leq1$ on $T^*M$ and $\pi^*\theta$ is a global $1$-form on $T^*M$, $X_t$ is 
uniquely defined on $T^*M$ by the above equation for $0\leq t\leq1.$ Next, we prove that $X_t$ is Lipschitz continuous for all $t\in[0,1].$
Notice that the twisted cotangent bundle $(T^*M,\omega)$ is strongly geometrically bounded (see \cite[Proposition 4.1]{lu_weinstein_1998-1}), there exists $\alpha>0$ such that 
$$\alpha\left\lVert X_t\right\rVert^2_g\leq\omega_t(X_t,J_tX_t)=\pi^*\theta[J_tX_t] $$
where $g$ is a standard metric induced by a Riemannian metric on $M$ and $J_t$ is an almost complex structure on $T(T^*M)$ such that $\omega_t(J_t\cdot,\cdot)=g(\cdot,\cdot).$ Thus
$$\left\lVert X_t\right\rVert_g\leq\frac{1}{\alpha}\sup_{(q,p)\in T^*M}\left\lvert \pi^*\theta(q,p)\right\rvert=\frac{1}{\alpha}\sup_{q\in M}\left\lvert \theta(q)\right\rvert <V$$
for some constant $V>0.$
Then the flow $\phi^t$ exists globally on $T^*M$ for all $t\in[0,1].$
$\hfill\qedsymbol$

%\begin{rmk}
%For a general closed magnetic form $\sigma$ on $M$, and a twisted symplectic form $\omega =\omega_0-\pi^*\sigma$, 
%the $T^*M.$ In fact, there is a neighborhood $N$ of $\text{Fix}(R)$ in $T^*M$, 
%a neighborhood $N_0$ of the 0-section $M_0$ of $T^*M$, and a diffeomorphism
%$$\Phi: N_0(\subset T^*M)\rightarrow N(\subset T^*M)$$
%such that $\Phi^*\omega=\omega_0$ and $ \Phi R_0=R\Phi.$ See \cite[Theorem 2]{meyerHamiltonianSystemsDiscrete1981}.
%\end{rmk}
~\\

The above symplectomorphism $\Phi$ is the well known momentum shift map (see \cite{guillemin_symplectic_1984}). 
By simple calculation, we obtain that $\Phi(q,p)=(q,p-\theta(q))$ in the local coordinates $(q_1,\dots,q_N)$.
We denote by $R_1:=\Phi\circ R_0\circ\Phi$, then $R_1^*\omega=-\omega$ is an anti-symplectic involution on $(T^*M,\omega)$.
Besides, we have $R_1(q,p)=\Phi\circ R_0\circ\Phi(q,p)=(q,-p-2\theta(q)),\forall (q,p)\in T^*M$.

%\subsection{The conjugation between the Hamiltonian flow \texorpdfstring{$\Psi_H^t$}{} on \texorpdfstring{$(T^*M,\omega, R)$}{} and the Lagrangian flow \texorpdfstring{$\Psi_{L_{\theta}}^t$}{} on \texorpdfstring{$TM$}{}}

Let $H: S_1\times T^*M\rightarrow\mathbb{R}$ be a 1-periodic Tonelli Hamilition.
Define $H_{\theta}=H\circ \Phi,$ then $H_{\theta}$ is also a 1-periodic Tonelli Hamilition. 
We denote by $X_K$ the Hamilition vector field of some Hamiltonian $K$ and $\Psi_K^t$ its Hamilition flow. 
Since $\Phi$ is symplectic , the Hamilition flow $\Psi_H^t$ of $H$ on $(T^*M,\omega)$ is conjugated to
the flow $\Psi_{H_{\theta}}^t$ of $H_{\theta}$ on $(T^*M,\omega_0)$  by the momentum shift map $\Phi$. 
Moreover, there is a one-to-one correspondence between the $\tau$-periodic brake orbits of $H$ in $(T^*M,\omega,R_1)$ 
with the $\tau$-periodic brake orbits of $H_{\theta}$ in $(T^*M,\omega_0,R_0)$.

As has been mentioned in Section \ref{sec2:2}, the Hamiltonian flow $\Psi_{H_{\theta}}^t$ is also conjugated to the Euler-Lagrangian flow $\Psi_{L_{\theta}}^t$
by Legendre transform. There is a one-to-one correspondence between the $\tau$-periodic brake orbits of $H_{\theta}$ in $(T^*M,\omega_0,R_0)$ 
and the $\tau$-periodic 
symmetric orbits of $L_{\theta}$ in the configuration space $M$. Therefore Theorem \ref{sec1:th3} is equivalent to Theorem \ref{sec1:th1} 
if $L$ is exactly the dual Lagrangian of $H$ and $\theta$ a primitive for $\sigma.$

Denote by $L, L_{\theta}: S_1\times TM\rightarrow\mathbb{R}$ the dual Tonelli Lagrangian of $H$ and $H_{\theta}$ respectively,
we have  
\begin{align*}
  L_{\theta}(t,q,v) & = \max_{p\in T_q^*M}\left\{p[v]-H_{\theta}(t,q,p)\right\} \\
  & = \max_ {p\in T_q^*M}\left\{p[v]-H(t,q,p-\theta(q))\right\}\\
  & = \max_ {p\in T_q^*M}\left\{(p+\theta(q))[v]-H(t,q,p)\right\}\\
  & = \max_ {p\in T_q^*M}\left\{p[v]-H(t,q,p)+\theta(q)[v]\right\}\\
  & = L(t,q,v)+\theta(q)[v],
  \end{align*}
and $L_{\theta}(-t,q,-v)=L_{\theta}(t,q,v).$

\section{Convex quadratic modification}

%The convex quadratic modification technique was introduced in \cite{abbondandoloHighActionOrbits2007}, where the authors showed how the multiplicity results holding for Lagrangian 
%which behave quadratically in $v$ for $v$ large can be entended to the classical setting of Tonelli Lagrangians.
We follow the mathod of \cite{mazzucchelli_lagrangian_2011} to modify the Tonelli Lagrangian
to be a convex quadratic-growth Lagrangian. 

Recall that $L_{\theta}=L+\theta$ is the dual Lagrangian of the Tonelli Hamilition $H_{\theta}$ on
$(T^*M,\omega_{0},R_0)$. 
Define $K:=\max_{q\in M}\sup_{|v|_q=1}|\theta(q)[v]|$, and
\begin{align*}
  C(H,\theta) & :=\max\{H(t,q,p)\ |\ (t,q,p)\in S_1\times T^*M,\, |p|_q\leq K+1\}\\
              & \geq\max\{H(t,q,p-\theta(q))\ |\ (t,q,p)\in S_1\times T^*M,\,|p|_q\leq1\}\\
              & =\max\{H_{\theta}(t,q,p)\ |\ (t,q,p)\in S_1\times T^*M,\,|p|_q\leq1\},
\end{align*}
we have
\begin{align*}
  L_{\theta}(t,q,v) & \geq \max_{|p|_q\leq 1}\{p[v]-H_{\theta}(t,q,p)\}\\
                    & \geq \max_{|p|_q\leq 1}\{p[v]\}-\max_{|p|_q\leq 1}\{H_{\theta}(t,q,p)\}\\
                    &\geq|v|_q-\max\{H_{\theta}(t',q',p')\ |\ 
(t',q',p')\in S_1\times T^*M, |p'|_{q'}\leq1\}\\
                    & \geq|v|_q-C(H,\theta).
\end{align*}
\begin{defe}[Convex quadratic modification] The Lagrangian
$L_T: S_1\times TM\rightarrow \mathbb{R}$ is a convex quadratic $T$-modification of a
Tonelli Lagrangian $L$ when

(M1)$L_T(t,q,v)=L(t,q,v)$ for each $(t,q,v)\in S_1\times TM$ with $|v|_q\leq T$;

(M2) $L_T$ is a convex quadratic-growth Lagrangian, i.e. satisfying \eqref{sec2:311};

(M3)$L_T(t,q,v)\geq|v|_q-C$ for each $(t,q,v)\in S_1\times TM$, and some constant $C\in\mathbb{R}$.
\end{defe}

For each real $T>0$, the convex quadratic $T$-modification of $L_\theta$, denoted by
$L_{\theta,T}: S_1\times TM\rightarrow \mathbb{R}$, exists and is also time-reversible, i.e.
$$L_{\theta,T}(-t,q,-v)=L_{\theta,T}(t,q,v),\ \forall t\in S_1.$$ 

For the reader's convenience, we present a brief introduction. The details are included in \cite{abbondandolo_high_2007}. 
Let $\phi: \mathbb{R}\rightarrow\mathbb{R}$ be a smooth increasing function such that $\phi(s)=s$ for $s\leq 1$ and $\phi(s)$ is constant
 for $s\geq 2$. Let $\lambda$ be a positive number such that 
 \begin{equation*}
   \lambda\geq \max\{L_\theta(t,q,v)\,|\, (t,q,v)\in S_1\times TM,\, |v|_q\leq 2T\},
 \end{equation*}
and $L_{1,T}:S_1\times TM\rightarrow\mathbb{R}$ be $L_{1,T}:=\lambda\phi(\tfrac{L_\theta}{\lambda}).$ Since $L_\theta$
satisfies (L1) in Section \ref{sec2:2} and $L_{1,T}$ is constant outside a compact set, we can find a positive number $\mu$ such that $\mu I+\partial_{vv}L_{1,T}>0$ on $S_1\times TM.$
We define $\psi:\mathbb{R}\rightarrow\mathbb{R}$ by $\psi(s)=0$ for $s\leq T^2$ and $\psi(s)=\mu s-2\mu T^2$ for $s\geq 4T^2$. 
Then the convex quadratic-growth $T$-modification can be constructed as
\begin{equation*}
  L_{\theta,T}:=L_{1,T}(t,q,v)+\psi(|v|_q^2).
\end{equation*}
In fact, if we choose $\mu$ satisfies 
\begin{equation*}
  4T\mu\geq1,\quad 2T^2\mu\geq 2T-C(H,\theta)-\min L_{1,T},
\end{equation*}
$L_{\theta,T}$ will satisfies (M1)(M2)(M3) (see \cite[Section 5]{abbondandolo_high_2007}).
The associated mean action functional is
$$EA_T(\xi)=\int_0^{1}L_{\theta,T}(t,\xi(t),\dot{\xi}(t))\text{d}t,\ \forall\xi\in EH(1).$$

Now we introduce the properties of $EA_T$ of the convex quadratic $T$-modification $L_{\theta,T}$.

\begin{prop}\label{sec5:11}
   Let $\gamma: S_{1}\rightarrow M$ be a symmetric $1$-periodic solution of the Euler-Lagrange system of $L_{\theta}$ and 
$U>\max\{|\dot{\gamma}(t)|_{\gamma(t)}\,|\, t\in S_{1}\}$, then

(1) $\gamma$ is a critical point of $EA_T$ for each $T\geq U$; 

(2)\begin{equation*}
  (m^-(EA_T,\gamma),m^0(EA_T,\gamma))=(m^-(EA_U,\gamma),m^0(EA_U,\gamma));
\end{equation*}

\end{prop}

\Proof. (1) We notice that  
$L_\theta(t,\gamma,\dot{\gamma})=L_{\theta,T}(t,\gamma,\dot{\gamma})$, then $\gamma$ is also a solution of the Euler-Lagrange system of $L_{\theta,T}$.
This implies that $\gamma$ is an extremal of $EA_T.$ 

\begin{rmk}
  Conversely, if $\gamma$ is a critical point of $EA_U$, $\gamma$ is in fact smooth and is a solution
  of the Euler-Lagrange system of $L_{\theta,U}$ (see e.g. \cite{fathi_weak_2008,buttazzo_one-dimensional_1998}).
   Furthermore, if we have $U>\max\{|\dot{\gamma}(t)|_{\gamma(t)} |\ t\in S_1\}$,
$\gamma$ will be a smooth symmetric $1$-periodic solution of the Euler-Lagrange system of $L_{\theta}$.
\end{rmk}

(2) Since $L_{\theta}$ and $L_{\theta,T}$ coincide along the 
lifted curve $(\gamma,\dot{\gamma}): S_{1}\rightarrow TM$ for each $T\geq U$, we have 
\begin{align}
  &\textnormal{Hess}EA_T(\gamma)[\xi,\eta] \notag\\
  =&\int_0^1[\partial_{vv}L_{\theta,T}(t,\gamma(t),\dot{\gamma}(t))\dot{\xi}\cdot\dot{\eta}
  +\partial_{qv}L_{\theta,T}(t,\gamma(t),\dot{\gamma}(t))\xi\cdot\dot{\eta}\notag\\
  &\qquad +\partial_{vq}L_{\theta,T}(t,\gamma(t),\dot{\gamma}(t))\dot{\xi}\cdot\eta
  +\partial_{qq}L_{\theta,T}(t,\gamma(t),\dot{\gamma}(t))\xi\cdot\eta]dt \notag\\
  =&\int_0^1[\partial_{vv}L_{\theta}(t,\gamma(t),\dot{\gamma}(t))\dot{\xi}\cdot\dot{\eta}
  +\partial_{qv}L_{\theta}(t,\gamma(t),\dot{\gamma}(t))\xi\cdot\dot{\eta} \notag\\
  &\qquad+\partial_{vq}L_{\theta}(t,\gamma(t),\dot{\gamma}(t))\dot{\xi}\cdot\eta
  +\partial_{qq}L_{\theta}(t,\gamma(t),\dot{\gamma}(t))\xi\cdot\eta]dt \label{sec4:12}
\end{align}
Therefore the Hessian of $EA_T$ at $\gamma$ is independent of the choice of $T$ for $T\geq U.$ In particular, we have
$$\textnormal{Hess}EA_U(\gamma)=\textnormal{Hess}EA_T(\gamma)$$
This implies the conclusion.$\hfill \qedsymbol$

~\\
There is a priori estimate for the symmetric critical point of $EA_T$, which is an immediate consequence of \cite[Lemma 5.2]{abbondandolo_high_2007}.

\begin{prop} \label{sec5:120}
 For each $\alpha>0$ and $m\in\mathbb{N}$, there exists $\tilde{T}=\tilde{T}(\alpha,m)>0$ such that, for any
 $T$-modification $L_{\theta,T}$ of $L_{\theta}$ with $T>\tilde{T}$ and for any $n \in\{1,\dots,m\}$,
 the following holds: if $\gamma$ is a critical point of $EA_T^{[n]}$ such that $EA_T^{[n]}\leq
 \alpha$, then $\max\{|\dot{\gamma}(t)|_{\gamma(t)}| t\in S_{n}\}\leq \tilde{T}$. In particular, $\gamma$ is a periodic 
 solution of the Euler-Lagrange system of $L_\theta$, and $EA^{[n]}(\gamma)=EA^{[n]}_T(\gamma).$
 $\hfill \qedsymbol$
\end{prop}

 ~\\
Let $\gamma^{[m]}=\psi^{[m]}(\gamma)$ be an isolated critical point of $EA^{[m]}_T$ with critical value $c=EA^{[m]}_T(\gamma^{[m]})$
and $\max\{|\dot{\gamma}(t)|_{\gamma(t)}\,|\, t\in S_1\}<U/2.$ In the following, we will consider the local homology group 
of $EA_T^{[m]}$ at $\gamma^{[m]}$ given by 
\begin{equation*}
  C_*(EA_T^{[m]},\gamma):=H_*((EA_T^{[m]})_c\cup\{\gamma^{[m]}\},(EA^{[m]}_T)_c),
\end{equation*}
 where $(EA_T^{[m]})_c:= (EA_T^{[m]})^{-1}(-\infty,c)$ is the sublevel of $EA_T^{[m]}$ at $c$, and $H_*$ is the singular homology factor with an arbitrary coefficient group.

Since the local homology group is of a local nature, we consider the 
local coordinates near $\gamma$ as in Section \ref{sec2:3}.

We give some notations first.
 Let $\rho$ be the injection radius of $M$, $B_\rho\subset \mathbb{R}^N$ be an open neighborhood of $0\in\mathbb{R}^N$
 with radius $\rho$. We denote by $V_\rho(m):=EW^{1,2}(S_m,B_\rho)$ and $EW(m):=EW^{1,2}(S_m,\mathbb{R}^N)$.
 Let
 $$X(m):=\{\text{ time-reversible loops in } C^1(S_m,\mathbb{R}^N)\}$$
 be the Banach space with the $C^1$-norm.
  Let $V_\rho^X(m):=V_\rho(m)\cap X(m)$.

 Since $\gamma$ is smooth, the local coordinates near $\gamma^{[m]}$ is given by  
\begin{equation}
  \Theta_{\gamma^{[m]}}: V_\rho(m)\rightarrow \mathcal{U}_{\gamma^{[m]}}\subset EW(m),\ \xi\mapsto \exp_{\gamma^{[m]}}(\Phi^{[m]}\xi),
\end{equation}
then \begin{align*}
  &EA^{[m]}_T\circ\Theta_{\gamma^{[m]}}(\xi)\\
 % =&\frac{1}{k1}\int_0^{k1} L\left(t,\Theta_\gamma(\xi)(t),\frac{\text{d}}{\text{d}t}\Theta(\gamma)(t)\right)\text{d}t\\
 =&\frac{1}{m}\int_0^{m} L_\theta\left(t,\Theta_{\gamma^{[m]}}(\xi)(t),\text{d}(\Theta_{\gamma^{[m]}})(\xi)\cdot\dot{\xi}+\frac{\text{d}}{\text{d}t}\Theta_{\gamma^{[m]}}(t)\circ\xi\right)\text{d}t\\
 =&\frac{1}{m}\int_0^{m}L_\theta\circ\pi_{\gamma^{[m]}}(t,\xi(t),\dot{\xi}(t))\text{d}t,
      \end{align*}
      
  We further choose $\rho$ small enough such that 
  \begin{equation*}
    \max_{t\in S_1,q\in\bar{B}_\rho} \left|\frac{d}{dt}\Theta_{\gamma}(q)(t)  \right| <U/2.
  \end{equation*}
Let 
\begin{equation}
  C_\rho:=\max_{t\in S_1,q\in\bar{B}_\rho}\sup_{|v|=1}\vert d\Theta_\gamma(q)(t)\cdot v\vert. \label{sec4:8}
\end{equation}  
Then for each $T\geq U$ and $\xi\in V_\rho(m)$ with $\text{ess sup}_{t\in S_m}|\dot{\xi}(t)|<T/2C_\rho=:\widehat{T}$, we have 
   \begin{equation*}
     \left\vert d\Theta_{\gamma^{[m]}}(\xi)\cdot \dot{\xi}(t)+\frac{d}{dt}\Theta_{\gamma^{[m]}}(t)\circ\xi  \right\vert<T.
   \end{equation*}
   Therefore,
 \begin{align*}
  L_{\theta,T}\circ\pi_{\gamma^{[m]}}(t,\xi(t),\dot{\xi}(t))&=L_{\theta}\circ\pi_{\gamma^{[m]}}(t,\xi(t),\dot{\xi}(t)),\\
  EA^{[m]}_T\circ\Theta_{\gamma^{[m]}}(\xi)&=EA^{[m]}\circ\Theta_{\gamma^{[m]}}(\xi).
\end{align*}
Here
$0=\Theta_\gamma^{-1}{\gamma}$ is the (smooth) critical point of both $EA^{[m]}\circ\Theta_{\gamma^{[m]}}$ and $EA^{[m]}_T\circ\Theta_{\gamma^{[m]}}$.
We denote by
 $$\widehat{L}=L_\theta\circ\pi_{\gamma^{[m]}},\,\widehat{L}_{\widehat{T}}=L_{\theta,T}\circ\pi_{\gamma^{[m]}}: S_m\times B_\rho\times\mathbb{R}^N\rightarrow \mathbb{R},$$ 
 and 
 $$\widehat{EA}^{[m]}_{\widehat{T}}=EA^{[m]}_T\circ\Theta_{\gamma^{[m]}}: EW(m)\rightarrow\mathbb{R}.$$ Then 
 \begin{equation}
   \widehat{L}_{\widehat{T}}(t,q,v)=\widehat{L}(t,{q},v),\quad \forall (t,q,v)\in S_1\times B_\rho\times\mathbb{R}^N\text{ with }|v|<\widehat{T}.\label{sec4:9}
 \end{equation}
Note that to simplify the notation, we simply write $L,L_T$ and $EA^{[m]}_T$ for $\widehat{L},\widehat{L}_{\widehat{T}}$ and $\widehat{EA}^{[m]}_{\widehat{T}}$ respectively.

~\\
Since $L_T$ is a convex quadratic-growth Lagrangian, $EA^{[m]}_T$ is $C^{2-0}$ on $EW(m)$, it is well known that  
$EA^{[m]}_T$ satisfies the generalized splitting lemma condition (S)(F)(B) (see in Appendix) given by Lu in \cite{lu_splitting_2013}.
%By the generalized shifting theorem (see in A??), the local homology is also finitely generated. 
%By the factor property of the splitting lemma, the iteration map $\psi^{[m]}$ induces the local homology isomorphism 
%if the Morse index and nullity of $\gamma$ are unchanged under iteration.
More precisely, we have the following conclusions.

Without loss of generality, we consider the case $m=1$,
 and denote by $V_\rho=V_\rho(1), X=X(1), V_\rho^X=V_\rho^X(1)$ and $EW=EW(1)$.

(S) Since $X$ is dense in $EW$ and 
$V_\rho^X$ is dense in $V_\rho.$ $X$ satisfies the condition (S).

(F1) Since $EA_T: V_\rho\rightarrow \mathbb{R}$ is a $C^{2-0}$ functional, $EA_T$ satisfies the condition (F1).

(F2) We denote by $EA_{T,X}:=EA_T|_{V_\rho^X}$. $EA_{T,X}$ is a $C^2$ functional. Let $\nabla EA_{T}\in EW$
 be the gradient of $EA_{T}$ with respect to the metric $\Angles{\cdot}{\cdot}$, i.e.
\begin{equation*}
  dEA_{T}(\gamma)[\xi]=dEA_{T,X}(\gamma)[\xi]=\Angles{\nabla EA_T(\gamma)}{\xi}\quad \forall \gamma\in V_\rho^X,\, \xi\in X.
\end{equation*}
Then the image of $\nabla EA_{T}$ is contained in $X$, and $\mathcal{A}_T:=\nabla EA_{T}|_{V_\rho^X}:V_\rho^X\rightarrow X$
is continuous differentiable (see \cite[Lemma 3.2]{lu_corrigendum_2009}).

 (F3) By \cite[Proposition 2.1]{jiang_generalization_1999}, there exists a uniformly continuous map $\mathcal{B}_T:V_\rho^X\rightarrow L_s(EW)$
 such that $\mathcal{B}_T(\gamma)$ is a bounded linear self-adjoint operatr on $EW$ for each $\gamma\in V_\rho^X$, and 
 \begin{equation*}
   d^2EA_{T,X}(\gamma)[\xi,\eta]=\Angles{\mathcal{B}_T(\gamma)\xi}{\eta}.
 \end{equation*}

(B) These conditions are verified in \cite[Lemma 3.4, 3.5]{lu_corrigendum_2009}.

By \eqref{sec4:12}, $\mathcal{B}_T(0)$ is independent of the choice $T$ for $T\geq U$. We denote by $\mathcal{B}:=\mathcal{B}_T(0)$. Let $EW^0:=\text{Ker}\mathcal{B}$,
$EW^-$ (resp. $EW^+$) be the negative definite (resp. positive definite) subspace of $\mathcal{B}$ in $EW$. 
Conditions in (B) imply that both $EW^0$ and $EW^1$ are finitely dimensional subspaces contained in $X$ (see \cite[Proposition B.2]{lu_splitting_2013}).
 Let $X=EW^0\oplus  X^{\pm }$, where $X^\pm:=X\cap EW^\pm=X^-+X\cap EW^+.$ Let $B_\epsilon$ be a neighborhood of 0 with radius $\epsilon$ in $EW$.
 We define $B^*_{\epsilon}:=B_\epsilon\cap EW^*$ for $*=0,+,-,\pm$.

 \begin{lem}[Generalized splitting lemma]\label{sec4:11}
  Under the above assumptions (S)(F)(B), for each $T\geq U$, if $m^0(EA_T,0)>0$, there exists 
  $\epsilon_T,\bar{\epsilon}$ with $0<\epsilon_T\leq\bar{\epsilon}$, 
   a (unique) $C^1$ 
  map $h: B_{\bar{\epsilon}}^0\rightarrow X^\pm$ satisfying $h(0)=0$ and 
  \begin{equation*}
    (I-P^0)\mathcal{A}_T(z+h(z))=0,\quad z\in B_{\bar{\epsilon}}^0,
  \end{equation*}
  an open neighborhood $W_T\subset B_{\bar{\epsilon}}^0$ of 0 and an origin-preserving homeomorphism 
  \begin{equation*}
    \Phi_T: B_{{\epsilon}_T}^0\times (B_{{\epsilon}_T}^++B_{{\epsilon}_T}^-)\rightarrow W_T
  \end{equation*}
  of form $\Phi_T(z,u^++u^-)=z+h(z)+\phi_{z,T}(u^++u^-)$ with $\phi_z(u^++u^-)\in H^\pm$ 
  with $ \phi_{z,T}(u^++u^-)\in EW^\pm $ such that 
  such that
  \begin{equation*}
    EA_T\circ\Phi_T(z,u^++u^-)=\Vert u^+\Vert^2-\Vert u^-\Vert^2+EA_T{(z+h(z))}
  \end{equation*}
  for all $(z,u^++u^-)\in B_{{\epsilon}_T}^0\times (B_{{\epsilon}_T}^++B_{{\epsilon}_T}^-),$ and that 
  \begin{equation*}
  \Phi_T(B_{{\epsilon}_T}^0)\times(B_{{\epsilon}_T}^++B_{{\epsilon}_T}^-)\subset X.
\end{equation*}
Moreover, the homeomorphism $\Phi_T$ also has properties:

(a1) For each $z\in B_{{\epsilon}_T}^0$, $\Phi_T(z,0)=z+h(z),\phi_{z,T}(u^++u^-)\in EW^-$ 
if and only if $u^+=0$.

(a2) The restriction of $\Phi_T$ to $B_{{\epsilon}_T}^0\times B_{{\epsilon}_T}^-$ is a homeomorphism from 
$B_{{\epsilon}_T}^0\times B_{{\epsilon}_T}^-\subset X\times X$ onto $\Phi_T(B_{{\epsilon}_T}^0\times B_{\epsilon_T}^-))\subset X$
even if the topologies on these two sets are chosen as the indued one by $X$.

The map $h$ and the function $B_{\bar{\epsilon}}^0\ni z\mapsto EA_T(z+h(z))$ also satisfy:

(b1) The map $h$ is $C^1$ in $V_\rho^X$, and
\begin{equation*}
  h'(z)=-[(I-P^0)\mathcal{A'}_T(z)|_{X^\pm}]^{-1}(I-P^0)\mathcal{A'}_T(z)\quad z\in B_{\bar{\epsilon}}^0;
\end{equation*}

(b2) $ \beta(z):=EA(z+h(z))=EA_T(z+h(z))$ for each $T\geq U$;

(b3) $\beta$ is $C^{2}$ for any $z\in B_{\bar{\epsilon}}^0$, and 
\begin{equation*}
  d\beta(z_0)(z)=\Angles{\mathcal{A}_T(z_0+h(z_0))}{z}\quad z_0\in B_{\bar{\epsilon}}^0;
\end{equation*}

(b4) If $\mathbf{0}$ is an isolated critical point of $EA_{T,X}$, then it is also an isolated critical point of $\beta$.
\end{lem}

\Proof. Since $EA_T$ is the mean action of a convex quadratic-growth Lagrangian $L_T$, by the generalized 
splitting lemma (see Appendix \ref{secA:1}), for each $T\geq U$, there exist 
 $\epsilon_T>0$, 
   a (unique) $C^1$ 
  map $h_T: B_{\epsilon_T}^0\rightarrow X^\pm$ satisfying $h_T(0)=0$ and 
  \begin{equation*}
    (I-P^0)\mathcal{A}_T(z+h_T(z))=0,\quad z\in B_{\epsilon_T}^0,
  \end{equation*}
  an open neighborhood $W_T\subset B_{\epsilon_T}^0$ of 0 and an origin-preserving homeomorphism 
  \begin{equation*}
    \Phi_T: B_{{\epsilon}_T}^0\times (B_{{\epsilon}_T}^++B_{{\epsilon}_T}^-)\rightarrow W_T
  \end{equation*}
  of form $\Phi_T(z,u^++u^-)=z+h_T(z)+\phi_{z,T}(u^++u^-)$
  satisfying the conditions (a1)(a2) and (b1)(b3)(b4).

  Here we only have to prove that $h_T$ is in fact independent of the choice $T.$ Let $\widehat{U}:=\frac{U}{2C_\rho}$, where $C_\rho$ is defined in \eqref{sec4:8}.
  We have the following 
  \begin{clm}
    There exists a positive $\bar{\epsilon}<\frac{\widehat{U}}{2}$, a unique map $h: B_{\bar{\epsilon}}^0\rightarrow X^\pm$
  such that 
  
  (i) $h(0)=0$ and $(I-P^0)(\mathcal{A}_T(z+h(z)))=0$ for all $z\in B_{\bar{\epsilon}}^0$;
  
  (ii) $\Vert h'(z)\Vert_{X^\pm}\leq1,\forall z\in \bar{B}_{\bar{\epsilon}}^0$.
  \end{clm}
  
  ~\\
   For each $z\in \bar{B}_{\bar{\epsilon}}^0\subset X$, we have $\max\{|\dot{z}(t)|\,|\, t\in S_1\}<\tfrac{\widehat{U}}{2}$. 
  Condition (ii) implies that 
  \begin{equation*}
    \vert \dot{z}(t)+h'(z)\dot{z}(t)\vert<\widehat{U}.
  \end{equation*}
Therefore, by \eqref{sec4:9} we have $EA_T(z+h(z))=EA(z+h(z))$ for all $z\in \bar{B}_{\bar{\epsilon}}^0$ and each $T\geq U$.

The assertion of the claim is implicitly contained in the proof of \cite[Lemma 3.1]{lu_splitting_2013}. 
Now we give the details. 
Under the assumptions in (B), it was proved in \cite{jiang_generalization_1999} that 
$\mathcal{B}(X^\pm)\subset X^\pm$ and $\mathcal{B}|_{X^\pm}:X^\pm\rightarrow X^\pm$
is an isomorphism. 

Let $B_{X,r}\subset V_\rho^X$ be a neighboehood of $\mathbf{0}\in X$ with radius $r$. 
If $0<r<\widehat{U}$, then for each $x\in B_{X,r}$, we have $\max\{|\dot{x}(t)|\,|\,t\in S_1\}<r<\widehat{U}\leq \widehat{T}$ and 
\begin{equation*}
  L_T(t,x(t),\dot{x}(t))=L(t,x(t),\dot{x}(t)).
\end{equation*}
Therefore,
\begin{equation*}
  \Vert \mathcal{A}_T(x_1)-\mathcal{B}x_1-\mathcal{A}_T(x_2)+\mathcal{B}x_2\Vert_X=\Vert \mathcal{A}_U(x_1)-\mathcal{B}x_1-\mathcal{A}_U(x_2)+\mathcal{B}x_2\Vert_X.
\end{equation*}
Since $\mathcal{A}_U$ is $C^1$, it is strictly F-differentiable at $\mathbf{0}\in X$. It follows that
\begin{equation*}
  \Vert \mathcal{A}_U(x_1)-\mathcal{B}x_1-\mathcal{A}_U(x_2)+\mathcal{B}x_2\Vert_X\leq K_r\Vert x_1-x_2\Vert_X
\end{equation*}
for all $x_1,x_2\in B_{X,r}$ with constant $K_r\rightarrow 0$ as $r\rightarrow 0$ (see the proof of [26,cor.3]).
Let $C_1=\Vert \mathcal{B}|_{X^\pm}\Vert^{-1}_{L(X^\pm,X^\pm)}$, $C_2=\Vert I-P^0\Vert_{L(X^\pm,X^\pm)}$.
We fix a small $0<\bar{\epsilon}<\frac{\widehat{U}}{2}$ so that $C_1C_2K_r<\frac{1}{2}$.
   
We consider the map 
\begin{equation*}
  S_T: B_{\bar{\epsilon}}\times(B_{X,\bar{\epsilon}}\cap X^\pm)\rightarrow X^\pm
\end{equation*}
given by 
\begin{equation*}
  S_T(z,x)=-(\mathcal{B}|_{X^\pm})^{-1}(I-P^0)\mathcal{A}_T(z+x)+x.
\end{equation*}
Let $z_1,z_2\in B_{\bar{\epsilon}}^0$ and $x_1,x_2\in B_{X,\bar{\epsilon}}\cap X^\pm$. Since $\mathcal{B}x_i\in X^\pm$
and $\mathcal{B}z_i=0,\, i=1,2,$ we get
\begin{align}
  &\Vert S_T(z_1,x_1)-S_T(z_2,x_2)\Vert_{X^\pm}\notag\\ 
  \leq &C_1C_2\Vert \mathcal{A}_T(z_1+x_1)-\mathcal{B}(z_1+x_1)-\mathcal{A}_T(z_2+x_2)+\mathcal{B}(z_2+x_2)\Vert_X            \notag\\ 
  =&C_1C_2\Vert \mathcal{A}_U(z_1+x_1)-\mathcal{B}(z_1+x_1)-\mathcal{A}_U(z_2+x_2)+\mathcal{B}(z_2+x_2)\Vert_X              \notag\\
  \leq& C_1C_2K_{\bar{\epsilon}}\Vert z_1+x_1-z_2-x_2\Vert_X           \notag\\ 
  <&\frac{1}{2}\Vert z_1+x_1-z_2-x_2\Vert_X,\quad\text{ if } (z_1,x_2)\neq(z_2,x_2).\label{sec4:10}
\end{align}
In particular, for any $z\in B_{\bar{\epsilon}}$ and $x_1,x_2\in B_{X,\bar{\epsilon}}\cap X^\pm$, it holds that 
\begin{equation*}
  \Vert S_T(z,x_1)-S_T(z,x_2)\Vert_{X^\pm}<\frac{1}{2}\Vert x_1-x_2\Vert_X,\quad x_1\neq x_2,
\end{equation*}
and $\Vert S_T(z,0)\Vert_{X^\pm}<\frac{1}{2}\Vert z\Vert_X<\frac{\bar{\epsilon}}{2}$.
Then it is well known that there exists a unique map 
\begin{equation*}
  h_T:B_{\bar{\epsilon}}^0\rightarrow B_{X,\bar{\epsilon}}\cap X^\pm
\end{equation*}
such that $S_T(z,h_T(z))=h_T(z)$ or equivalently 
\begin{equation*}
  (I-P^0)\mathcal{A}_T(z+h_T(z))=0,\quad \forall z\in B_{\bar{\epsilon}}.
\end{equation*}
Clearly, $h_T(0)=0$. By the defnition of $S_T$ and \eqref{sec4:10}, we have 
\begin{equation*}
  \Vert h_T(z_1)-h_T(z_2)\Vert_X <\Vert z_1-z_2\Vert_X,\quad \forall z_1,z_2\in B_{\bar{\epsilon}}^0.
\end{equation*}
Since $h_T$ is differentiable on $B_{\bar{\epsilon}}^0$, we have 
\begin{equation*}
  \Vert h_T'(z)\Vert_{X^\pm}\leq1,\forall z\in \bar{B}_{\bar{\epsilon}}^0,
\end{equation*}
which imples that $ \vert \dot{z}(t)+h_T'(z)\dot{z}(t)\vert<\widehat{U},$ and
$$L_T(t,z(t)+h_T(z)(t),\dot{z}(t)+h_T'(z)\dot{z}(t))=L_U(t,z(t)+h_T(z)(t),\dot{z}(t)+h_T'(z)\dot{z}(t)).$$
$h_T$ satisfies the equation
\begin{equation*}
  (I-P^0)\mathcal{A}_U(z+h_T(z))=0,\quad \forall z\in B_{\bar{\epsilon}}^0
\end{equation*}
as well as $h_U$. By the uniqueness of the solution of the above equation, we must have $h_T=h_U$, 
$\hfill \qedsymbol$

\begin{cor}[shifting]\label{sec4:4}
  Let $\gamma$ be an isolated critical point of $EA_T$ with the property that $\max\{|\dot{\gamma}(t)|_{\gamma(t)}\,|\, t\in S_1\}<U/2.$
  For each $T\geq U$, we have 
  \begin{equation*}
    C_*(EA_T,\gamma;\mathbb{G})\cong C_{*-m^-(EA_U,0)}(\beta,0;\mathbb{G}).
  \end{equation*}
  where $\beta(z)=EA\circ\Theta_\gamma(z+h(z))$, $\mathbb{G}$ is an Abelian group. In particular, $C_*(EA_T,\gamma)=0$ if 
  \begin{equation}
    *<m^-(EA_U,\gamma)\text{ or }*>m^-(EA_U,\gamma)+m^0(EA_U,\gamma).\label{sec4:2}
  \end{equation}
\end{cor}
\Proof.
  By Lemma \ref{sec4:11} and the Shifting Theorem \ref{secA:3}, we have the homology isomorphism 
  \begin{equation*}
   C_*(EA_T,\gamma;\mathbb{G})\cong C_*(EA_T\circ\Theta_\gamma,0;\mathbb{G})\cong C_*(\alpha+\beta,0;\mathbb{G})\cong C_{*-m^-(EA_T,\gamma)}(\beta,0;\mathbb{G}).
  \end{equation*}
  Therefore, $C_*(EA_T,\gamma;\mathbb{G})=0$ if 
  \begin{equation*}
    *<m^-(EA_T,\gamma)\text{ or }*>m^-(EA_T,\gamma)+m^0(EA_T,\gamma).
  \end{equation*}
  \eqref{sec4:2} is obtained by Proposition \ref{sec5:11}(ii).
  $\hfill \qedsymbol$

  ~\\
   We denote by 
  \begin{equation}
    \mathcal{X}(m)=\{\text{ time-reversible loops in }C^1(S_m,M)\}\label{sec4:7}
  \end{equation}
  the $C^2$-Banach manifold with topology of uniformly convergence of time-reversible loops and their derivatives.
  Then $\mathcal{X}(m)\subset EH(m)$ is dense in $EH(m)$ and the inclusion $\mathcal{X}(m)\hookrightarrow EH(m)$ is continuous.
   $\gamma^{[m]}\in \mathcal{X}(m)$ if $\gamma$ is a solution of the Euler-Lagrange system of $L_\theta.$
  We will show that the inclusion map between  $\mathcal{X}(m)$ and $EH(m)$ 
induces an isomorphism between the local homology of $EA^{[m]}_T|_{\mathcal{X}(m)}$ 
and $EA^{[m]}_T$ at $\gamma$.
In particular, we show that the local homology of $EA^{[m]}_T$ at $\gamma$ is independent of the choice $T$ for $T\geq U$.

  We simply write $ \mathcal{X}= \mathcal{X}(1)$  and 
  \begin{equation*}
    EA^{[m]}_{U,\mathcal{X}}:=EA^{[m]}_U|_{\mathcal{X}(m)}.
  \end{equation*}

  Let $\widehat{U}:=\frac{U}{2C_\rho}$, $\mathcal{W}(m):=\Theta_{\gamma^{[m]}}(B_{\widehat{U}}(m))$
and $\mathcal{U}(m):=\Theta_{\gamma^{[m]}}(B^X_{\widehat{U}}(m))$, where
$B_{\widehat{U}}(m)=\{\xi\in V_\rho(m)\,|\, \Vert\xi\Vert<\widehat{U}\}$ and
$B_{\widehat{U}}^X(m)=\{\xi\in V^X_\rho(m)\,|\, \Vert{\xi}\Vert_X<\widehat{U}\}$.
Since $\Vert\xi\Vert\leq\Vert \xi\Vert_X$, we have $\mathcal{U}(m)\subset\mathcal{W}(m)$ and $\mathcal{U}(m)\subset\mathcal{W}^ {\mathcal{X}}(m):=\mathcal{W}(m)\cap{\mathcal{X}}(m) $
  (with the topology induced as a subset of ${\mathcal{X}}(m)$)
   is an open subset of $\mathcal{W}^{ \mathcal{X}}(m)$. Therefore, we have 
   \begin{equation*}
     EA_{U, \mathcal{X}}^{[m]}|_{\mathcal{U}(m)}=EA_{T, \mathcal{X}}^{[m]}|_{\mathcal{U}(m)}=EA_{T, \mathcal{X}}^{[m]}|_{\mathcal{W}^{ \mathcal{X}}(m)}.
   \end{equation*}
   In particular, for $m=1$ we have the following two inclusion maps,
   \begin{align*}
     I_\mathcal{X}:& ((EA_{U, \mathcal{X}}|_\mathcal{U})_c\cup\{\gamma\},(EA_{U, \mathcal{X}}|_\mathcal{U}))
                   \rightarrow ((EA_{T, \mathcal{X}}|_{\mathcal{W}^{ \mathcal{X}}})_c\cup\{\gamma\}, (EA_{T, \mathcal{X}}|_{\mathcal{W}^{ \mathcal{X}}})_c),\\
     I_{\mathcal{X},EH}:& ((EA_{T, \mathcal{X}}|_{\mathcal{W}^{ \mathcal{X}}})_c\cup\{\gamma\}, (EA_{T, \mathcal{X}}|_{\mathcal{W}^{ \mathcal{X}}})_c)
                   \rightarrow ((EA_T|_{\mathcal{W}})_c\cup\{\gamma\}, (EA_{T}|_{\mathcal{W}})_c).
   \end{align*}
   By the excision property, $I_\mathcal{X}$ induces homology isomorphism
   \begin{equation}
    I_\mathcal{X*}: C_*(EA_{U, \mathcal{X}}|_\mathcal{U},\gamma)\xrightarrow{\simeq} C_*(EA_{T, \mathcal{X}}|_{\mathcal{W}^{ \mathcal{X}}},\gamma).\label{sec4:6}
   \end{equation}

  \begin{prop}\label{sec4:5}
   Let $\gamma$ be an isolated critical point of $EA_T$ with critical value $c:=EA_T(\gamma)$.
  Suppose $\max\{|\dot{\gamma}(t)|_{\gamma(t)}\,|\, t\in S_1\}<U/2$, then the inclusion map 
    \begin{equation*}
      J:=I_{\mathcal{X},EH}\circ I_\mathcal{X}: (((EA_{U,\mathcal{X}}|_\mathcal{U})_c\cup\{\gamma\}),(EA_{U,\mathcal{X}}|_\mathcal{U}))
      \rightarrow (((EA_T|_{\mathcal{W}})_c\cup\{\gamma\}), (EA_{T}|_{\mathcal{W}})_c)
    \end{equation*}
  induces the homology isomorphism with coefficients in a field $\mathbb{F}$
\begin{equation}
  J_*: C_*(EA_{U,\mathcal{X}},\gamma;\mathbb{F})\xrightarrow{\simeq} C_*(EA_{T},\gamma;\mathbb{F}).\label{sec4:3}
\end{equation}
  \end{prop}
\Proof. We have the following commutative diagram
\begin{center}
  \begin{tikzcd} 
    ((EA_{U, \mathcal{X}}\circ\Theta_\gamma|_{B^X_{\widehat{U}}})_c\cup\{0\},(EA_{U, \mathcal{X}}\circ\Theta_\gamma|_{B^X_{\widehat{U}}})_c)         \arrow[d,"\Theta_\gamma"]   \arrow[r,"I_{X,EW}"] 
  & ((EA_T\circ\Theta_\gamma|_{B_{\widehat{U}}})_c\cup\{0\}, (EA_{T}\circ\Theta_\gamma|_{B_{\widehat{U}}})_c)              \arrow[d,"\Theta_\gamma"]  \\
  ((EA_{U, \mathcal{X}}|_\mathcal{U})_c\cup\{\gamma\},(EA_{U, \mathcal{X}}|_\mathcal{U})_c)       \arrow[r,"I_{\mathcal{X},EH}"]      
  & ((EA_T|_{\mathcal{W}})_c\cup\{\gamma\}, (EA_{T}|_{\mathcal{W}})_c),  
  \end{tikzcd}
\end{center}
where $I_{X,EW}$ is the inclusion map induced by the continuous inclusion $X\hookrightarrow EW.$

Since $EA_{U, \mathcal{X}}\circ\Theta_\gamma$ satisfies the conditions (S)(F)(B), 
Theorem \ref{secA:2} implies that $I_{X,EW}$ induces an isomorphism in homology with coefficients in $\mathbb{F}$. 
It is easy to know that the coordinates map $\Theta_\gamma$ induces an isomorphism in homology. Therefore,
$I_{\mathcal{X},EH}$ also induces isomorphism in homology with coefficients in $\mathbb{F}$.
$J_*$ will be a composition of two homology isomorphisms $I_{X,EW*}$ and $ I_\mathcal{X*}$
$\hfill \qedsymbol$

\begin{prop}\label{sec4:1}
  Let $\gamma$ be an isolated critical point of $EA_T$ with critical value $c:=EA_T(\gamma)$.
  Suppose $\gamma$ is also an isolated critical 
  point of $EA_T^{[m]}$ for some $m\in\mathbb{N}$. If 
  \begin{equation*}
    (m^-(EA_T,\gamma),m^0(EA_T,\gamma))=(m^-(EA_T^{[m]},\gamma^{[m]}),m^0(EA_T^{[m]},\gamma^{[m]})),
  \end{equation*}
  then the iteration map $\psi^{[m]}$, restricted as a map of pairs of the form
  $$\psi^{[m]}: ((EA_T)_c\cup\{\gamma\},(EA_T)_c)\hookrightarrow((EA_T^{[m]})_c\cup\{\gamma\},(EA_T^{[m]})_c),$$
  induces the homology isomorphism
  \begin{equation}
    \psi^{[m]}_*: C_*(EA_T,\gamma)\xrightarrow{\simeq}C_*(EA_T^{[m]},\gamma).\label{sec4:311}
        \end{equation}
\end{prop}
The statement is a consequence of the factor property of the generalized splitting lemma in \cite{lu_splitting_2013}.
We refer the reader to \cite{lu_conley_2009} for a proof.

\section{Proof of Theorem \ref{sec1:th3}}
  Let $\mathbb{K}:=\{2^n\, |\, n\in \mathbb{N}\cup\{0\} \}$ and $a\in \mathbb{R}$ be a constant greater than 
$$\max_{q\in M}\left\{\int_{0}^1L(t,q,0)dt\right\}.$$
We assume that symmetric solutions of the Euler-Lagrangian system of $L_\theta=L+\theta$ 
with period in $\mathbb{K}$ and mean action less than $a$ are finite, then we will deduce a contradiction. 
We denote these orbits by
  $$\gamma_1,\dots,\gamma_r.$$
 Let $2^n$ be the minimum of their commen periods in $\mathbb{K}$.  
 It is well known that $\tilde{\gamma}:\mathbb{R}\rightarrow M$ is 
 a $j$-periodic solution of the Euler-lagrangian system of the following time-rescaled Tonelli Lagrangian 
 \begin{equation*}
   \widetilde{L}_1(t,q,v):=L_1(2^nt,q,2^{-n}v)
 \end{equation*} 
 if and only if $\gamma:=\tilde{\gamma}(2^{-n}\cdot{}):\mathbb{R}\rightarrow M$ is a $2^nj$-periodic solution of the Euler-Lagrangian system of ${L_1}$.
 Moreover, $\tilde{\gamma}$ and $\gamma$ have the same mean action with respect to $\widetilde{L}_1$ and $L_1$ respectively.
Without loss of generality, we assume in the following that all these orbits $\gamma_1,\dots,\gamma_r$ have period 1.
Moreover, we consider the homology groups in this section that are assumed to have coefficients in the field $\mathbb{Z}_2.$
\begin{lem} For each $T>0$ and $n\in\mathbb{N}$, the homology of the sublevel $(EA_T^{[n]})_a$
is non-trivial in degree $N$, i.e.
$$H_N((EA_T^{[n]})_a)\neq0.$$
\end{lem}

\Proof. By our choice of the constant $a$ in Theorem \eqref{sec1:th3}, we have 
\begin{equation*}
  EA^{[n]}_T\circ \text{pt}^{[n]}(q)=EA^{[n]}\circ \text{pt}^{[n]}(q)<a,\quad \forall q\in M,
\end{equation*}
where $\text{pt}^{[n]}: M\hookrightarrow EH(n)$ is the embedding that maps a point to the constant loop at that point.
Then $\text{pt}^{[n]}(q)\in (EA^{[n]})_a.$
%We denote by 
%\begin{equation*}
%  \text{ev}: EH(n)\rightarrow M,\quad \gamma\mapsto \gamma(0).
%\end{equation*}
%the evaluation map. Then $\text{ev}\circ \text{pt}^{[n]}=id_M$. 
It follows that the homology homomorphism 
\begin{equation*}
  \text{pt}^{[n]}_*: H_N(M)\rightarrow H_N((EA_T^{[n]})_a)
\end{equation*}
induced by $\text{pt}^{[n]}$ is a monomorphism. Since $H_N(M)\neq0$, the Lemma follows.
$\hfill \qedsymbol$

~\\
Let $U/2>\max\{|\dot{\gamma}_\nu(t)|_{\gamma(t)}\ |\ t\in S_1,\nu=1,\dots,r\}$. 
Then for each $m\in\mathbb{K}$, the Morse index and nullity of $EA_T^{[m]}$ at each $\gamma_\nu$
is independent of the choice $T$ for $T\geq U.$
We denote by 
\begin{equation*}
  (m^-(EA^{[m]},\gamma^{[m]}_\nu),m^0(EA^{[m]},\gamma^{[m]}_\nu)):=(m^-(EA^{[m]}_T,\gamma^{[m]}_\nu),m^0(EA^{[m]}_T,\gamma^{[m]}_\nu)).
\end{equation*}
Let $\{\gamma_1,\dots,\gamma_s\}$ be the symmetric orbits with mean Morse index zero, 
and $\{\gamma_{s+1},\dots,\gamma_r\}$ be the symmetric orbits with positive mean index. We choose $m_1\in\mathbb{K}$
large enough such that 
\begin{equation*}
  m^-(EA^{[m_1]},\gamma_\nu)\geq N+1,\ \forall \nu\in\{s+1,\dots,r\}.
\end{equation*}
Then by Corollary \ref{sec4:4}, we have
\begin{equation*}
  C_N(EA_{T}^{[m_1]},\gamma_\nu)=0,\ \forall \nu\in \{s+1,\dots,r\} 
  \text{ and }T\geq U.
\end{equation*}
By the Morse inequality, we have 
\begin{equation*}
  0\neq H_N((EA^{[m_1]}_T)_a)\leq\sum_{\nu=1}^sC_N(EA^{[m_1]}_T,\gamma_\nu)
\end{equation*}
thus $\exists \gamma\in\{\gamma_1,\cdots,\gamma_s\}$ such that $C_N(EA^{[m_1]}_T,\gamma)\neq0$.
Since $\widehat{m}^-(EA,\gamma)=0$,  \eqref{sec2:33} implies that 
\begin{equation*}
  0\leq m^-(EA^{[m]},\gamma)+m^0(EA^{[m]},\gamma)\leq N,\ \forall m\in\mathbb{K}.
\end{equation*}
We choose an infinite subset $\mathbb{K}'\subset\mathbb{K}$ such that for each $m\in\mathbb{K}'$,
\begin{equation*}
 (m^-(EA^{[m]},\gamma),m^0(EA^{[m]},\gamma)= m^-(EA^{[m_1]},\gamma),m^0(EA^{[m_1]},\gamma)).
\end{equation*}
Up to a time-rescaling of the Lagrange function $L_\theta$, we assume $1\in\mathbb{K'}$.

%%%%%%%%

By Proposition \ref{sec4:1}, the iteration map $\psi^{[m]}$ induces homology isomorphism 
\begin{equation}
  \psi_*^{[m{}]}: 0\neq C_N(EA_{T}^{},\gamma)\xrightarrow{\simeq} C_N(EA_{T}^{[m]},\gamma),\ \forall T\geq U
\end{equation}

Let $c:=EA(\gamma)=EA_{T}^{}(\gamma^{})=EA_{T}^{[m]}(\gamma^{[m]})$, and let $\epsilon>0$ be small enough so that $c+\epsilon<a$ and there is no
$\gamma_\nu\in\{\gamma_1,\dots,\gamma_r\}$ such that $EA(\gamma_\nu)\in(c,c+\epsilon).$

%By Theorem \ref{sec5:120}, since $T_1>\tilde{T}(a,m_1),$ we know that the only critical points of $EA_{T_1}^{}$ in the open sublevel 
%$(EA_{T_1}^{})_a$ are $\gamma_1^{},\dots,\gamma_r^{}.$ However, it may not be true for the sublevel $EA_{T_1}^{[m]}$ 

%It is noted that if $T<\tilde{T}(a,m),$ there may be
%other critical points except $\gamma_1^{[m]},\dots,\gamma_r^{[m]}$ in $(EA_{T}^{[m]})_a$. 
%So the inclusion
%$$((EA_{T}^{[m]})_{c}\cup\{\gamma^{[m]}\}, (EA_{T}^{[m]})_{c})\hookrightarrow((EA_{T}^{[m]})_{c+\epsilon},(EA_{T}^{[m]})_{c})$$
%may do not induce a monomorphism in homology.

By Proposition \ref{sec5:120}, if we choose $m\in\mathbb{K}'$ and $T>\max\{U,\tilde{T}(a,m)\},$ 
there is no critical point of $EA_{T}^{[m]}$ with  
 value in $(c, c+\epsilon)$. This implies that the inclusion 
\begin{equation}
  ((EA_T^{[m]})_{c}\cup\{\gamma^{[m]}\}, (EA_T^{[m]})_{c})\hookrightarrow((EA_T^{[m]})_{c+\epsilon},(EA_T^{[m]})_{c})\label{sec6:3}
\end{equation}
 induces a monomorphism in homology (See \cite[Chapter 1, Theorem 4.2]{chang_infinite_1993}).
 
 By Proposition \ref{sec4:5}, the inclusion map
 \begin{equation}
 	\begin{aligned}
 		J^{[m]}: 
 		(((EA^{[m]}_{U,\mathcal{X}}|_{\mathcal{U}(m)})_c\cup\{\gamma^{[m]}\}),(EA^{[m]}
 		_{U,\mathcal{X}}|_{\mathcal{U}(m)}))\rightarrow\\ 
 		(((EA^{[m]}_T|_{\mathcal{W}(m)})_c\cup\{\gamma^{[m]}\}), 
 		(EA^{[m]}_{T}|_{\mathcal{W}(m)})_c)\label{sec6:2}
 	\end{aligned}
 \end{equation}
induces the homology isomorphism.

By \eqref{sec6:3}\eqref{sec6:2}, we obtain a monomorphism in homology
\begin{equation}\label{sec6:6}
  \rho^{[m]}_*: C_*(EA^{[m]}_{U,\mathcal{X}}|_{\mathcal{U}(m)},\gamma)\rightarrow H_*((EA_{T}^{[m]})_{c+\epsilon},(EA_{T}^{[m]})_{c})
\end{equation}
for each $m\in\mathbb{K}'$ and $T>\max\{U,\tilde{T}(a,m)\}$.

Notice that the iteration map $\psi^{m{}}: EH(m_1)\rightarrow EH(m)$ 
restricts to an inclusion map on every action sublevel $(EA_T)_c$ of $EA_T$ and 
commutes with $\rho^{[m]}$. Namely, for each $m\in\mathbb{K}'$ and 
$T>\max\{U,\tilde{T}(a,m)\}$,

\begin{center}
  \begin{tikzcd}
    ((EA_{U,\mathcal{X}}^{})_{c}\cup\{\gamma^{}\}, (EA_{U,\mathcal{X}}^{})_{c})     \arrow[d,"\rho^{}"]       \arrow[r,"\psi^{[m{}]}"]      
  & ((EA_{U,\mathcal{X}}^{[m]})_{c}\cup\{\gamma^{[m]}\}, (EA_{U,\mathcal{X}}^{[m]})_{c})    \arrow[d,"\rho^{[m]}"]  \\
  ((EA_T^{})_{c+\epsilon},(EA_T^{})_{c})     \arrow[r,"\psi^{[m{}]}"]       
  & ((EA_T^{[m]})_{c+\epsilon},(EA_T^{[m]})_{c}).               
  \end{tikzcd}
\end{center}
 This
induces the following commutative diagram in homology
\begin{center}
  \begin{tikzcd} 
     C_*(EA_{U,\mathcal{X}}^{},\gamma)              \arrow[d,"\rho^{}_*"] \arrow[r,"\psi_*^{[m{}]}"] 
  & C_*(EA_{U,\mathcal{X}}^{[m]},\gamma)              \arrow[d,"\rho^{[m]}_*"]  \\
    H_*((EA_{T}^{})_{c+\epsilon},(EA_{T}^{})_{c})        \arrow[r,"\psi_*^{[m{}]}"]      
  & H_*((EA_{T}^{[m]})_{c+\epsilon},(EA_{T}^{[m]})_{c}).  
  \end{tikzcd}
\end{center}

 We will prove
\begin{lem}[Homological vanishing]\label{sec6:7}
  Let $[\mu]\in C_*(EA_{U,\mathcal{X}}^{},\gamma) $. Then
   there exists $\bar{m}=\bar{m}(L_{\theta},[\mu])\in \mathbb{K'}$ and $\bar{T}=\bar{T}(L_\theta,[\mu])$ such that
   for each $T\geq \bar{T}$,
  \begin{equation*}
    \psi_*^{[\bar{m}{}]}\circ\rho_*^{}[\mu]=0 \text{ in } H_*((EA_{T}^{[\bar{m}]})_{c+\epsilon},(EA_{T}^{[\bar{m}]})_{c}).
  \end{equation*}
\end{lem}

Theorem \ref{sec1:th3} readily follows from this lemma. Indeed, 
let $0\neq [\mu]\in C_N(EA_{U,\mathcal{X}}^{},\gamma) $. 
%Since $N>0$, $\gamma^{}$ can not be a local minimum of $\beta_{U}^{}$.
 We take $m\in\{m\in\mathbb{K'}\,|\, m\geq\bar{m}\}$ and 
$T\geq\max\{U,\tilde{T}(a,m),\bar{T}\}$, thus 
\begin{equation*}
  \psi_*^{[m{}]}\circ\rho^{}[\mu]=\psi_*^{[m/\bar{m}]}\circ\left(\psi_*^{[\bar{m}{}]}\circ\rho_*^{}[\mu]\right)=0
\end{equation*}
in $H_N((EA_{T}^{[m]})_{c+\epsilon},(EA_{T}^{[m]})_{c})$.
However, the above commutative diagram implies that 
\begin{equation*}
  (0\neq)\rho_*^{[m]}\circ\psi_*^{[m{}]}[\mu]=\psi_*^{[m{}]}\circ\rho_*^{}[\mu],
\end{equation*}
which is a contradiction.
$\hfill \qedsymbol$

\section{Proof of Lemma \ref{sec6:7}}\label{sec7:5}
   
The proof of Lemma \ref{sec6:7} is based on a homotopical technique that is essentially
due to Bangert. It is noted that the Bangert homotopy we need here should be time-reversible 
which is different from that in 
\cite{bangert_homology_1983,long_multiple_2000,lu_conley_2009,mazzucchelli_lagrangian_2011}. 

Let $L: S_1\times TM\rightarrow\mathbb{R}$ be a convex quadratic-growth Lagrangian that satisfies
\begin{equation*}
  L(-t,q,-v)=L(t,q,v),\quad (t,q,v)\in S_1\times TM.
\end{equation*}
$EA: EH(1)\rightarrow\mathbb{R}$ is the mean action functional of $L$. 
\begin{lem}[Time-reversible Bangert homotopy] \label{sec7:3}
  Let $c_1<c_2\leq\infty$ and $\sigma: \Delta^q\rightarrow EH(1)$
be a continuous singular simplex such that 
\begin{align*}
  \sigma(\Delta^q)\subset (EA)_{c_2},\quad
  \sigma(\partial\Delta^q)\subset (EA)_{c_1}.
\end{align*}
Then there exists 
$\bar{n}=\bar{n}(L,\sigma)\in\mathbb{N}$ and, for every integer 
$n\geq\bar{n}$, a homotopy 
\begin{equation*}
  B_{\sigma}^{[2n]}: [0,1]\times(\Delta^q,\partial\Delta^q)\rightarrow ((EA)_{c_2},(EA)_{c_1}),
\end{equation*}
which we call Bangert homotopy, such that 

(i) $B_{\sigma}^{[2n]}(0,\cdot)=\sigma^{[2n]}:=\psi^{[2n]}\circ\sigma$,

(ii) $B_{\sigma}^{[2n]}(1,\Delta^q)\subset(EA^{[2n]})_{c_1},$

(iii) $B_{\sigma}^{[2n]}(\cdot,z)=B_{\sigma}^{[2n]}(0,z)$ if $z\in\partial\Delta^q.$
\end{lem}
%\begin{rmk}
  
% With the time-reversible Bangert homotopy Lemma \ref{sec7:3}, we have the general homological vanishing theorem.
 
% Let $c_1<c_2\leq\infty,$
%   where the sublevel $(EA)_{c_1}$ is not empty,   and let $[\mu]\in\text{H}_*((EA)_{c_2},(EA)_{c_1}).$ Then
%   there exists $\bar{m}=\bar{m}(L,[\mu])\in \mathbb{K}$such that $\psi_*^{[\bar{m}]}[\mu]=0$ in$\text{H}_*((EA^{[\bar{m}]})_{c_2},(EA^{[\bar{m}]})_{c_1}).$
%\end{rmk}
~\\
\Proof. First of all, we introduce some notations. Let $\alpha_i: [a_i,b_i]\rightarrow M$ for $i=1,2$ be two paths defined on bounded intervals 
with $\alpha_1(b_1)=\alpha_2(a_2).$ We define the inverse path $\alpha_1^{-1}$ as a path going along  $\alpha_1$ from $\alpha_1(b_1)$ to $\alpha_1(a_1).$
We denote by $\overline{\alpha_1}:
[a,b]\rightarrow M$ the reparametrization of $\alpha_1$ on $[a,b]$ with the formula 
$$\overline{\alpha_1}(t):=\alpha\left(\tfrac{t-a}{b-a}\cdot(b_1-a_1)+a_1\right).$$

  Besides, $\alpha_1*\alpha_2$ is the concatenation path defined by 
first traversing $\alpha_1$ and then traversing $\alpha _2$.
And a parametrization of $\alpha_1*\alpha_2$ on $[a_1,b_1+b_2-a_2]$
%, with $\bar{\alpha}_2$ defined on $[b_1,b_1+b_2-a_2]$, is given by
% $\alpha _1$ can be parametrized as a map  $\alpha_1: [a_1,b_1]\rightarrow M$, 
is given by 
$$\alpha_1*\alpha_2(t):=
\begin{cases}
  \alpha_1(t),&t\in[a_1,b_1).\\
  \overline{\alpha_2}(t),&t\in[b_1,b_1+b_2-a_2].
\end{cases}$$

Now we construct a proper homotopy in $EH(2n).$ Consider a continuous map $\theta: [x_0,x_1]\rightarrow EH(1),$
where $[x_0,x_1]\subset\mathbb{R}$. For each $n\in\mathbb{N},$ we define $\theta^{[2n]}:=\psi^{[2n]}\circ\theta: [x_0,x_1]\rightarrow EH(2n).$
We will build another continuous map $\theta^{\left\langle 2n\right\rangle}$ which is a crucial step to construct the Bangert homotopy. 
 
Denote by $\text{ev}: W^{1,2}(S_1,M)\rightarrow M$ the evaluation map, given by 
$$\text{ev}(\xi)=\xi(0),\quad \forall\xi\in W^{1,2}(S_1,M).$$
This map is smooth, which implies that the initial point curves $\text{ev}\circ\theta: [x_0,x_1]\rightarrow M$
is $\text{uniformly continuous}$. 

%$\textbf{injectivity radius}$ In particular, there exists a contant $\rho=\rho(\theta)\geq0$ such that,
%for each $x, x'\in[x_0,x_1]$ with $|x-x'|\leq\rho,$ we have that $\text{dist}(\text{ev}\circ\theta(x),\text{ev}\circ\theta(x'))$
%is less than the injectivity radius of $M$. Here, the "dist" is the Riemannian distance on $M$.

\textbf{Step 1.}\ If $\text{dist}(\text{ev}\circ\theta(x),\text{ev}\circ\theta(x'))$ is less than the injecivity radius, we define $\theta_x^{x'}$ 
as the shortest geodesic that connects the points $\text{ev}\circ\theta(x)$ and $\text{ev}\circ\theta(x').$ Notice that this geodesic
depends smoothly on its endpoints. For each $x\in [x_0,x_1],\;\text{dist}(\text{ev}\circ\theta(x_0),\text{ev}\circ\theta(x))$ may be 
greater than the injectivity radius. Here we choose $l\in\{0,1,\dots,l_1\}$ such that 
$x_0+l\rho\leq x\leq x_0+(l+1)\rho$ where $l_1$ satisfies $x_0+l_1\rho\leq x_1\leq x_0+(l_1+1)\rho$. Then we can define
 the broken geodesics $\theta_{x_0}^{x}: [x_0,x]\rightarrow M$ 
and $\theta_x^{x_1}:[x,x_1]\rightarrow M$ by
\begin{align*}
  \theta_{x_0}^{x}&:=\theta_{x_0}^{x_0+\rho}*\theta_{x_0+\rho}^{x_0+2\rho}*\dots*\theta_{x_0+l\rho}^x,\\
  \theta_x^{x_1}&:=\theta_x^{x_0+(l+1)\rho}*\theta_{x_0+(j+1)\rho}^{x_0+(j+2)\rho}*\dots*\theta_{x_0+l_1\rho}^{x_1}.
\end{align*}

\textbf{Step 2.}\ For each $x\in [x_0,x_1]$, we define a loop $\tilde{\theta}^{\left\langle 2n\right\rangle }(x)$ in the following way. 
Let $l\in\{1,2,\dots,n-2\}$ and $y\in [0,\frac{x_1-x_0}{n}]$, we put 
\begin{align*}
  \tilde{\theta}^{\left\langle 2n\right\rangle}\!(x_0\!+\!y)&:=\theta^{[n-1]}(x_0)*\theta_{x_0}^{x_0+ny}*\theta^{[2]}(x_0\!+\!ny)\\
  &\quad *(\theta_{x_0}^{x_0+ny})^{-1}*\theta^{[n-1]}(x_0),\\
  \tilde{\theta}^{\left\langle 2n\right\rangle}\!\left(x_0\!+\!\tfrac{l}{n}(x_1\!-\!x_0)\!+\!y\right)&:=\theta^{[n-l-1]}(x_0)*\theta_{x_0}^{x_0+ny}*\theta(x_0\!+\!ny)*\theta_{x_0+ny}^{x_1}*\theta^{[2l]}(x_1)\\
  &\quad *(\theta_{x_0+ny}^{x_1})^{-1}*\theta(x_0\!+\!ny)*(\theta_{x_0}^{x_0+ny})^{-1}*\theta^{[n-l-1]}(x_0),\\
  \tilde{\theta}^{\left\langle 2n\right\rangle}\!\left(x_0\!+\!\tfrac{n-1}{n}(x_1\!-\!x_0)\!+\!y\right)&:=\theta(x_0\!+\!ny)*\theta_{x_0+ny}^{x_1}*\theta^{[2n-2]}(x_1)\\
  &\quad *(\theta_{x_0+ny}^{x_1})^{-1}*\theta(x_0\!+\!ny).
\end{align*}
Note that the loops $$\tilde{\theta}^{\left\langle 2n\right\rangle}\!(x_0\!+\!y),\ \tilde{\theta}^{\left\langle 2n\right\rangle}\!\left(x_0\!+\!\tfrac{l}{n}(x_1\!-\!x_0)\!+\!y\right),\
 \tilde{\theta}^{\left\langle 2n\right\rangle}\!\left(x_0\!+\!\tfrac{n-1}{n}(x_1\!-\!x_0)\!+\!y\right)$$ are defined respectively on 
$$[0,2n(1+y)],\ [0,2n+2(x_1-x_0)],\ [0,2n+2(x_1-x_0-ny)].$$

\textbf{Step 3.}
For each $y\in [0,\frac{x_1-x_0}{n})$ and $x=x_0+\tfrac{l}{n}(x_1-x_0)+y\in[x_0,x_1),$ we reparametrize the loop $\tilde{\theta}^{\langle 2n \rangle}(x)$ as a 
time-reversible one $\theta^{\left\langle 2n\right\rangle}(x)$ in $EH(2n).$ More precisely,
%Here, we take $y\in[0,\frac{x_1-x_0}{n}).$
\begin{equation*}
  \theta^{\left\langle 2n\right\rangle}\!(x_0\!+\!y)|_{[0,n]}:=
\begin{cases}
  \theta^{[n-1]}(x_0),& t\in[0,n-1]\\
  \overline{\theta_{x_0}^{x_0+ny}},& t\in[n-1,n-1+\frac{ny}{ny+1}]\\
  \overline{\theta(x_0\!+\!ny)},& t\in[n-1+\frac{ny}{ny+1},n].
\end{cases}
\end{equation*}

\begin{equation*}
  \theta^{\left\langle 2n\right\rangle}\!\left(x_0\!+\!\tfrac{1}{n}(x_1\!-\!x_0)\!+\!y\right)\!|_{[0,n]}:=
\begin{cases}
  \theta^{[n-2]}(x_0),& t\in[0,n-2]\\
  \overline{\theta_{x_0}^{x_0+ny}},& t\in[n-2,n-2+\frac{ny}{ny+1}]\\
  \overline{\theta(x_0\!+\!ny)},& t\in[n-2+\frac{ny}{ny+1},n-1]\\
  \overline{\theta_{x_0+ny}^{x_1}}, &t\in[n-1,n-1+\frac{x_1-x_0-ny}{x_1-x_0-ny+1}]\\
  \overline{\theta(x_1)},&t\in[n-1+\frac{x_1-x_0-ny}{x_1-x_0-ny+1},n]
\end{cases}
\end{equation*}

\begin{equation*}
  \theta^{\left\langle 2n\right\rangle}\!\left(x_0\!+\!\tfrac{l}{n}(x_1\!-\!x_0)\!+\!y\right)\!|_{[0,n]}:=
\begin{cases}
  \theta^{[n-l-1]}(x_0),& t\in[0,n-l-1]\\
  \overline{\theta_{x_0}^{x_0+ny}},& t\in[n-l-1,n-l-1+\frac{ny}{ny+1}]\\
  \overline{\theta(x_0\!+\!ny)},& t\in[n-l-1+\frac{ny}{ny+1},n-l]\\
  \overline{\theta_{x_0+ny}^{x_1}}, &t\in[n-l,n-l+\frac{x_1-x_0-ny}{x_1-x_0-ny+1}]\\
  \overline{\theta(x_1)},&t\in[n-l+\frac{x_1-x_0-ny}{x_1-x_0-ny+1},n-l+1]\\
  \overline{\theta^{[l-1]}(x_1)},&t\in [n-l+1,n].
\end{cases}
\end{equation*}
for $l=2,\dots,n-2.$

\begin{equation*}
  \theta^{\left\langle 2n\right\rangle}\!\left(x_0\!+\!\tfrac{n-1}{n}(x_1\!-\!x_0)\!+\!y\right)\!|_{[0,n]}:=
\begin{cases}
  \theta(x_0\!+\!ny),& t\in[0,1]\\
  \overline{\theta_{x_0+ny}^{x_1}}, &t\in[1,1+\frac{x_1-x_0-ny}{x_1-x_0-ny+1}]\\
  \overline{\theta(x_1)},&t\in[1+\frac{x_1-x_0-ny}{x_1-x_0-ny+1},2]\\
  \overline{\theta^{[n-2]}(x_1)},&t\in [2,n].
\end{cases}
\end{equation*}
and $\theta^{\left\langle 2n\right\rangle}(x_1)|_{[0,n]}:=\theta^{[2n]}(x_1).$ For $t\in [-n,0],$ we define
$$\theta^{\left\langle 2n\right\rangle}(x)(t):=\theta^{\left\langle 2n\right\rangle}(x)(-t),\quad \forall x\in[x_0,x_1].$$
Since $\theta^{\left\langle 2n\right\rangle}(x)(n)=\theta^{\left\langle 2n\right\rangle}(x)(-n)=0,\ 
\theta^{\left\langle 2n\right\rangle}(x) $ 
 can be extended as a loop with period $2n$. 
 It is easy to know that $\theta^{\left\langle 2n\right\rangle}(x)$ is a continuous loop on $S_{2n}$ and
  is a reparametrization of $\tilde{\theta}^{\left\langle 2n\right\rangle}(x)$
 for each $x\in[x_0,x_1]$. (See the example $\theta^{\left\langle 
 6\right\rangle}(x)$ in Figure \ref{fig1}).

 For fixed $l\in\{-n+1,-n+2,\dots,n-1\},$ $\theta^{\left\langle 2n\right\rangle}\left(x_0+\tfrac{l}{n}(x_1-x_0)+y\right)$ is continuous with respect to $y\in[0,\frac{x_1-x_0}{n}).$
 Besides, 
 $$\lim_{y\rightarrow(x_1-x_0)/n}\theta^{\left\langle 2n\right\rangle}\left(x_0+\tfrac{l}{n}(x_1-x_0)+ny\right)=\theta^{\left\langle 2n\right\rangle}\left(x_0+\tfrac{l+1}{n}(x_1-x_0)\right)$$
 for each $l\in\{-n+1,-n+2,\dots,n-1\}.$ Thus, $\theta^{\left\langle 2n\right\rangle}: [x_0,x_1]\rightarrow EH(2n)$ is a continuous map. 
\begin{figure}
	\begin{tikzpicture}[scale=0.6]
		%%椭圆 两个半轴长 1cm 和 1.5cm, 间距0.03cm, 拐弯处角度差20度
		\begin{scope}[xshift=1cm]
		\draw [wbl] (1,0)++(7*0.03,0)++({0.03+1-cos(10)},{1.5*sin(10)}) 
			coordinate (a0) arc (170:-170:1 and 1.5) to [out=90+10,in=0-30] 
			++({-0.03-1+cos(10)-0.075},{1.5*sin(10)}) to [out=180-30, in=-90+5] 
			++({-0.03-1+cos(10)-0.075},{1.5*sin(10)}) arc (10:350:1 and 1.5) to 
			[out=90-5, in=0-60] ++({-0.03-1+cos(10)},{1.5*sin(10)}) to 
			[out=180-60, 
			in=-90+5] ++({-0.03-1+cos(10)},{1.5*sin(10)}) coordinate (a1);
			\draw [wbl] (a1) arc (10:350:1 and 1.5) to [out=90-5, in=0-60] 
			++({-0.03-1+cos(10)},{1.5*sin(10)}) to [out=180-60, in=-90+5] 
			++({-0.03-1+cos(10)},{1.5*sin(10)}) coordinate (a2);
			\draw [wbl] (a0) to [out=-90-5, in=0+60] 
			++({-0.03-1+cos(10)},{-1.5*sin(10)}) to [out=180+60, in=90+5] 
			++({-0.03-1+cos(10)},{-1.5*sin(10)}) arc (-170:170:1 and 1.5)  to 
			[out=-90-5, in=0+60] ++({-0.03-1+cos(10)},{-1.5*sin(10)}) to 
			[out=180+60, in=90+5] ++({-0.03-1+cos(10)},{-1.5*sin(10)}) arc 
			(-170:-100:1 and 1.5);
			\draw [wbl] (a2) arc (10:350:1 and 1.5) to [out=90-5, in=180+30] 
			++({0.03+1-cos(10)+0.075},{1.5*sin(10)}) to [out=0+30, in=-90-5] 
			++({0.03+1-cos(10)+0.075},{1.5*sin(10)}) arc (170:-100:1 and 1.5);
			\node at (-3.5,0) {$\theta^{\langle6\rangle}(0)$};
		\end{scope}

		\begin{scope}[yshift=-4cm,xshift=1cm]
			\draw [wbl] 
			(1,0)++(5*0.03,0)++({0.03+1-cos(10)},{1.5*sin(10)})++(3,0) 
			coordinate (b0) arc (170:-170:1 and 1.5) to [out=90+5,in=0-60] 
			++({-0.03-1+cos(10)},{1.5*sin(10)}) to [out=180-60, in=-90+5] 
			++({-0.03-1+cos(10)},{1.5*sin(10)}) arc (10:350:1 and 1.5) to 
			[out=90-5, in=0] ++({-0.03-1+cos(10)},{1.5*sin(10)-0.03}) 
			coordinate (b1);
			\draw [wbl] (b0)  to [out=-90-10, in=0] 
			++({-0.03-1+cos(10)-0.03*2},{-1.5*sin(10)+0.03}) -- ++(-3+0.03*5,0) 
			coordinate (b3);
			\draw [wbl] (b1) -- ++(-3,0) to [out=180,in=-90+5] 
			++({-0.03-1+cos(10)-0.03*2},{1.5*sin(10)+0.03}) arc (10:350:1 and 
			1.5) to [out=90-5, in=0-60] ++({-0.03-1+cos(10)},{1.5*sin(10)}) to 
			[out=180-60, in=-90+5] ++({-0.03-1+cos(10)},{1.5*sin(10)}) 
			coordinate (b2);
			\draw [wbl] (b3)  to [out=180, in=90-5] 
			++({-0.03-1+cos(10)-0.03*2},{-1.5*sin(10)-0.03}) arc (-170:170:1 
			and 1.5) to [out=-90-5, in=60] ++({-0.03-1+cos(10)},{-1.5*sin(10)}) 
			to [out=180+60, in=-90+5] ++({-0.03-1+cos(10)},{-1.5*sin(10)}) arc 
			(-170:-100:1 and 1.5);
			\draw [wbl] (b2) arc (10:350:1 and 1.5) to [out=90-5, in=180+60] 
			++({0.03+1-cos(10)+0.03*2},{1.5*sin(10)}) to [out=60, in=-90+5] 
			++({0.03+1-cos(10)+0.03},{1.5*sin(10)}) arc (170:-100:1 and 1.5);
			\node at (-3.5,0)
			{$\theta^{\langle6\rangle}(y)|_{(0,\frac13)}$};
%			\node [below=1ex] at (-3.5,0) {$y \in ()$};
			
			\fill 
			(1,0)++(6*0.03,0)++(3,0)++(1,0)++({cos(120)},{1.5*sin(120)})circle 
			(2pt) node [above=1ex] {\footnotesize 
			$\theta^{\langle6\rangle}(y)(-t)$} ;
			\fill 
			(1,0)++(4*0.03,0)++(3,0)++(-1,0)++({cos(-60)},{1.5*sin(-60)})circle 
			(2pt) node [below=1ex] {\footnotesize 
			$\theta^{\langle6\rangle}(y)(t)$};
			 \draw [-stealth] 
			 (1,0)++(6*0.03,0)++(3,0)++(1,0)++({cos(60)},{1.5*sin(60)})++(0.3,0)
			  -- ++(2em,0);
			\draw [-stealth] 
			(1,0)++(6*0.03,0)++(3,0)++(1,0)++({cos(60)},{1.5*sin(60)})++(0.45,-0.5)
			 -- ++(2em,0);
			  \draw [-stealth] 
			 (1,0)++(6*0.03,0)++(3,0)++(1,0)++({cos(60)},{1.5*sin(60)})++(0.55,-1)
			 -- ++(2em,0);
			  \draw [-stealth] 
			 (1,0)++(6*0.03,0)++(3,0)++(1,0)++({cos(60)},{1.5*sin(60)})++(0.55,-1.5)
			 -- ++(2em,0);
			   \draw [-stealth] 
			 (1,0)++(6*0.03,0)++(3,0)++(1,0)++({cos(60)},{1.5*sin(60)})++(0.45,-2)
			 -- ++(2em,0);
			   \draw [-stealth] 
			 (1,0)++(6*0.03,0)++(3,0)++(1,0)++({cos(60)},{1.5*sin(60)})++(0.3,-2.5)
			 -- ++(2em,0);
		\end{scope}
		
		\begin{scope}[yshift=-8cm,xshift=1cm]
			\draw [wbl] 
			(1,0)++(5*0.03,0)++({0.03+1-cos(10)},{1.5*sin(10)})++(9,0) 
			coordinate (b0) arc (170:-170:1 and 1.5) to [out=90+5,in=0-60] 
			++({-0.03-1+cos(10)},{1.5*sin(10)}) to [out=180-60, in=-90+5] 
			++({-0.03-1+cos(10)},{1.5*sin(10)}) arc (10:350:1 and 1.5) to 
			[out=90-5, in=0] ++({-0.03-1+cos(10)},{1.5*sin(10)-0.03}) 
			coordinate (b1);
			\draw [wbl] (b0)  to [out=-90-10, in=0] 
			++({-0.03-1+cos(10)-0.03*2},{-1.5*sin(10)+0.03}) -- ++(-9+0.03*5,0) 
			coordinate (b3);
			\draw [wbl] (b1) -- ++(-9,0) to [out=180,in=-90+5] 
			++({-0.03-1+cos(10)-0.03*2},{1.5*sin(10)+0.03}) arc (10:350:1 and 
			1.5) to [out=90-5, in=0-60] ++({-0.03-1+cos(10)},{1.5*sin(10)}) to 
			[out=180-60, in=-90+5] ++({-0.03-1+cos(10)},{1.5*sin(10)}) 
			coordinate (b2);
			\draw [wbl] (b3)  to [out=180, in=90-5] 
			++({-0.03-1+cos(10)-0.03*2},{-1.5*sin(10)-0.03}) arc (-170:170:1 
			and 1.5) to [out=-90-5, in=60] ++({-0.03-1+cos(10)},{-1.5*sin(10)}) 
			to [out=180+60, in=-90+5] ++({-0.03-1+cos(10)},{-1.5*sin(10)}) arc 
			(-170:-100:1 and 1.5);
			\draw [wbl] (b2) arc (10:350:1 and 1.5) to [out=90-5, in=180+60] 
			++({0.03+1-cos(10)+0.03*2},{1.5*sin(10)}) to [out=60, in=-90+5] 
			++({0.03+1-cos(10)+0.03},{1.5*sin(10)}) arc (170:-100:1 and 1.5);
			\node at (-3.5,0) {$\theta^{\langle6\rangle}(\frac13)$};
		\end{scope}
		
		\begin{scope}[yshift=-12cm,xshift=1cm]
			\draw [wbl] 
			(1,0)++(4*0.03,0)++({0.03+1-cos(10)},{1.5*sin(10)})++(9,0) 
			coordinate (c0) arc (170:-170:1 and 1.5) to [out=90, in=-60] 
			++({-0.03-1+cos(10)},{1.5*sin(10)}) to [out=180-60, in=-90] 
			++({-0.03-1+cos(10)},{1.5*sin(10)}) arc (10:350:1 and 1.5) to 
			[out=90-5, in=0] ++({-0.03-1+cos(10)},{1.5*sin(10)-0.03}) 
			coordinate (c1);
			\draw [wbl] (c0) to [out=-90-5,in=0] 
			++({-0.03-1+cos(10)},{-1.5*sin(10)+0.03}) -- ++(-6+0.03*3,0) to 
			[out=180,in=90-5] ++({-0.03-1+cos(10)-0.03*3},{-1.5*sin(10)-0.03}) 
			arc (190:270:1 and 1.5) coordinate (c2);
			\draw [wbl] (c1) -- ++(-6+0.03*3,0) to [out=180,in=-90] 
			++({-0.03-1+cos(10)-0.03*3},{1.5*sin(10)+0.03}) arc (10:350:1 and 
			1.5) to [out=90-5, in=0] ++({-0.03-1+cos(10)},{1.5*sin(10)-0.03}) 
			coordinate (c3);
			\draw [wbl] (c3) -- ++(-3+0.03*3,0) to [out=180,in=-90] 
			++({-0.03-1+cos(10)-0.03*3},{1.5*sin(10)+0.03}) arc (10:350:1 and 
			1.5);
			\draw [wbl] (c2) arc (-90:170:1 and 1.5) to [out=-90-5,in=0] 
			++({-0.03-1+cos(10)},{-1.5*sin(10)+0.03}) -- ++(-3+0.03*3,0) to 
			[out=180,in=90-5] ++({-0.03-1+cos(10)-0.03*3},{-1.5*sin(10)-0.03}) 
			arc (190:270:1 and 1.5) coordinate (c4);
			\draw [wbl] (c4) arc (-90:170:1 and 1.5) to [out=-90-5,in=60] 
			++({-0.03-1+cos(10)},{-1.5*sin(10)}) to [out=180+60,in=90-5] 
			++({-0.03-1+cos(10)},{-1.5*sin(10)});
			\node at (-3.5,0) 
			{$\theta^{\langle6\rangle}(\cdot)|_{(\frac13,\frac23)}$};
%			\node [below=1ex] at (-3.5,0) {$y \in $};
			%node[此处设置node的类型]
			 \draw [-stealth] 
			(1,0)++(6*0.03,0)++(3,0)++(1,0)++({cos(60)},{1.5*sin(60)})++(0.25,0)
			-- ++(2em,0);
			\draw [-stealth] %[thick]
			(1,0)++(6*0.03,0)++(3,0)++(1,0)++({cos(60)},{1.5*sin(60)})++(0.4,-0.5)
			-- ++(2em,0);
			\draw [-stealth] 
			(1,0)++(6*0.03,0)++(3,0)++(1,0)++({cos(60)},{1.5*sin(60)})++(0.5,-1)
			-- ++(2em,0);
			\draw [-stealth] 
			(1,0)++(6*0.03,0)++(3,0)++(1,0)++({cos(60)},{1.5*sin(60)})++(0.5,-1.5)
			-- ++(2em,0);
			\draw [-stealth] 
			(1,0)++(6*0.03,0)++(3,0)++(1,0)++({cos(60)},{1.5*sin(60)})++(0.4,-2)
			-- ++(2em,0);
			\draw [-stealth] 
			(1,0)++(6*0.03,0)++(3,0)++(1,0)++({cos(60)},{1.5*sin(60)})++(0.25,-2.5)
			-- ++(2em,0);
		\end{scope}
		
		\begin{scope}[yshift=-16cm,xshift=1cm]
			\draw [wbl] 
			(1,0)++(-0.03*3,0)++({0.03+1-cos(10)},{1.5*sin(10)})++(9,0) 
			coordinate (d0) to [out=-90-5, in=180-60] 
			++({0.03+1-cos(10)},{-1.5*sin(10)}) to [out=-60, in=90+5] 
			++({0.03+1-cos(10)},{-1.5*sin(10)}) arc (-10:-350:1 and 1.5)  to 
			[out=-90-5, in=180-60] ++({0.03+1-cos(10)+0.03},{-1.5*sin(10)}) to 
			[out=-60, in=90+5] ++({0.03+1-cos(10)+0.03},{-1.5*sin(10)}) arc 
			(-170:170:1 and 1.5) coordinate (d1);
			\draw [wbl] (d0) arc (10:350:1 and 1.5) to [out=90-5, in=0] 
			++({-0.03-1+cos(10)+0.015-0.03*3},{1.5*sin(10)-0.03}) coordinate 
			(d2);
			\draw [wbl] (d2) -- ++(-9+0.03*6,0) to [out=180, in=-90+5] 
			++({-0.03-1+cos(10)+0.015-0.03*3},{1.5*sin(10)+0.03});
			\draw [wbl] (d1) to [out=-90+5, in=60] 
			++({-0.03-1+cos(10)},{-1.5*sin(10)}) to [out=180+60, in=90-5] 
			++({-0.03-1+cos(10)},{-1.5*sin(10)}) arc (-170:170:1 and 1.5) to 
			[out=-90+5, in=0] ++({-0.03-1+cos(10)-0.03*3},{-1.5*sin(10)+0.03}) 
			-- ++(-9+0.03*4,0) coordinate (d4);
			\draw [wbl] (d4) to [out=180, in=90+5] 
			++({-0.03-1+cos(10)-0.03*2},{-1.5*sin(10)-0.03}) arc (-170:170:1 
			and 1.5)  to [out=-90, in=60] ++({-0.03-1+cos(10)},{-1.5*sin(10)}) 
			to [out=180+60, in=90-5] ++({-0.03-1+cos(10)},{-1.5*sin(10)}) arc 
			(-10:-350:1 and 1.5);
			\node at (-3.5,0) {$\theta^{\langle6\rangle}(\frac23)$};
		\end{scope}
		\begin{scope}[yshift=-20cm,xshift=1cm]
			\draw [wbl] (10,0)++(7*0.03,0)++({0.03+1-cos(10)},{1.5*sin(10)}) 
			coordinate (a0) arc (170:-170:1 and 1.5) to [out=90+10,in=0-30] 
			++({-0.03-1+cos(10)-0.075},{1.5*sin(10)}) to [out=180-30, in=-90+5] 
			++({-0.03-1+cos(10)-0.075},{1.5*sin(10)}) arc (10:350:1 and 1.5) to 
			[out=90-5, in=0-60] ++({-0.03-1+cos(10)},{1.5*sin(10)}) to 
			[out=180-60, in=-90+5] ++({-0.03-1+cos(10)},{1.5*sin(10)}) 
			coordinate (a1);
			\draw [wbl] (a1) arc (10:350:1 and 1.5) to [out=90-5, in=0-60] 
			++({-0.03-1+cos(10)},{1.5*sin(10)}) to [out=180-60, in=-90+5] 
			++({-0.03-1+cos(10)},{1.5*sin(10)}) coordinate (a2);
			\draw [wbl] (a0) to [out=-90-5, in=0+60] 
			++({-0.03-1+cos(10)},{-1.5*sin(10)}) to [out=180+60, in=90+5] 
			++({-0.03-1+cos(10)},{-1.5*sin(10)}) arc (-170:170:1 and 1.5)  to 
			[out=-90-5, in=0+60] ++({-0.03-1+cos(10)},{-1.5*sin(10)}) to 
			[out=180+60, in=90+5] ++({-0.03-1+cos(10)},{-1.5*sin(10)}) arc 
			(-170:-100:1 and 1.5);
			\draw [wbl] (a2) arc (10:350:1 and 1.5) to [out=90-5, in=180+30] 
			++({0.03+1-cos(10)+0.075},{1.5*sin(10)}) to [out=0+30, in=-90-5] 
			++({0.03+1-cos(10)+0.075},{1.5*sin(10)}) arc (170:-100:1 and 1.5);
			\node at (-3.5,0) {$\theta^{\langle6\rangle}(1)$};
		\end{scope}	
	\end{tikzpicture}
	\caption{Description of $\theta^{\left\langle 6\right\rangle}(y)$}
\label{fig1}
\end{figure}

 We notice that the map $\tilde{\theta}^{\left\langle n\right\rangle}$ defined in \cite{mazzucchelli_lagrangian_2011} 
 can not be reparametrized as a symmetric loop, since the \emph{moving loop} $\theta(x_0+ny)$
 is 1-periodic, while a symmetric loop must go through the \emph{moving loop} $\theta(x_0+ny)$ twice at some time $t$ and $n-t$. 
 That is why we define a $2n$-periodic loop $\tilde{\theta}^{\left\langle 
 2n\right\rangle}$ and consider
  two \emph{moving loops} $\theta(x_0\!+\!ny)$ in $\tilde{\theta}^{\left\langle 2n\right\rangle}$.
  % $\theta^{\left\langle 2\right\rangle}(x_0+ny).$

For each $x\in[x_0,x_1]$, we define $\widehat{\theta}(x): S_2\rightarrow M$ as the loop obtained erasing from the formula of 
$\tilde{\theta}^{\left\langle 2n\right\rangle}$ the terms of powers of $\theta(x_0)$ and $\theta(x_1)$ and reparametrizing on $[0,2]$. Notice that $\widehat{\theta}$
is independent of the integer $n\in\mathbb{N}$ and, for each $x\in [x_0,x_1]$, the action
$EA(\widehat{\theta}(x))$ is finite and depends continuous on $x$. In particular we have
\begin{align*}
  C(\theta)&:=\max_{x\in[x_0,x_1]}\{EA(\widehat{\theta}(x))\}\\
           &\,=\max_{x\in[x_0,x_1]}\left\{\int_0^2L\left(t,\widehat{\theta}(x)(t),\tfrac{\text{d}}{\text{d}t}\widehat{\theta}(x)(t)\right)\text{d}t\right\}<\infty.
\end{align*}  
For each $n\in\mathbb{N},$ we have the estimate
\begin{equation}\label{sec7:6}
  \begin{aligned}
EA^{[2n]}(\theta^{\langle 2n\rangle}(x))&\leq\tfrac{1}{2n}[(n-2)\max\{EA(\theta(x_0)),EA(\theta(x_1))\}+EA(\widehat{\theta}(x))]\\
                         &\leq\max\{EA(\theta(x_0)),EA(\theta(x_1))\}+\tfrac{C(\theta)}{2n}.
\end{aligned}
\end{equation}

\textbf{Step 4.}
Now we construct the homotopy with respect to the $q$-singular simplex $\sigma:\Delta^q\rightarrow (EA)_{c_2}$.
Let $\mathbb{L}\subset\mathbb{R}^q$ be a straight line passing through the origin and the barycenter of $\Delta^q\subset\mathbb{R}^q.$ 
Then  $\mathbb{L}^\perp \oplus\mathbb{L}$ is a orthogonal decomposition of $\mathbb{R}^q.$ We can write $\mathbf{z}\in\Delta^q$
as $(\mathbf{y},x)\in\mathbb{L}^\perp \oplus\mathbb{L}.$ For each $s\in[0,1],$ let 
$s\Delta^q=\{s\mathbf{z}\ |\ \mathbf{z}\in\Delta^q\}$. For each fixed $\mathbf{y}\in\mathbb{L}^\perp$, we denote by 
$[x_0(\mathbf{y},s),x_1(\mathbf{y},s)]\subset\mathbb{L}$ the maximal interval such that $\{\mathbf{y}\}\times[x_0(\mathbf{y},s),x_1(\mathbf{y},s)]\subset s\Delta^q.$
Notice that 
$$\sigma(\mathbf{y},\cdot)|_{[x_0(\mathbf{y},s),x_1(\mathbf{y},s)]}: [x_0(\mathbf{y},s),x_1(\mathbf{y},s)]\rightarrow(EA)_{c_2}$$
 is a continuous map, we have the related 
time-reversible map $\sigma(\mathbf{y},\cdot)^{\langle 2n\rangle}$ for each $n\in\mathbb{N}.$ Now we define 
the Bangert homotopy $B_\sigma^{[2n]}(s,\mathbf{z}):[0,1]\times\Delta^q\rightarrow EH(2n)$ as
\begin{equation*}
  B_\sigma^{[2n]}(s,\mathbf{z}):=
  \begin{cases}
    {\sigma(\mathbf{y},)|_{[x_0(\mathbf{y},s),x_1(\mathbf{y},s)]}}^{\langle 2n\rangle}(x) &\mathbf{z}=(\mathbf{y},x)\in s\Delta^q,\\
    \sigma^{[2n]}(\mathbf{z}) & \mathbf{z}=(\mathbf{y},x)\in\Delta^q\setminus s\Delta^q.
  \end{cases}
\end{equation*}
Since $\sigma(\Delta^q)\subset(EA)_{c_2}$ and $\sigma(\partial\Delta^q)\subset(EA)_{c_1}$, by \eqref{sec7:6}
 there exists $\epsilon>0$ such that 
\begin{equation*}
  EA^{[2n]}(\sigma^{[2n]}(\Delta^q))<c_2-\epsilon,\ EA^{[2n]}(\sigma^{[2n]}(\partial\Delta^q))<c_1-\epsilon.
\end{equation*}
It is easy to know that $B_\sigma^{[2n]}$ is relative $\partial\Delta^q$,
we have 
\begin{equation*}
  EA(B_\sigma^{[2n]}(s,\partial\Delta^q))<c_1-\epsilon,\ \forall s\in[0,1].
\end{equation*}

Let\begin{equation*}
  C(\sigma):=\max\{C(\sigma(\mathbf{y},\cdot)|_{[x_0(\mathbf{y},s),x_1(\mathbf{y},s)]})\ |\ s\in[0,1],(\mathbf{y},x)\in\Delta^q\}.
\end{equation*}
It is easy to compute that 
\begin{equation*}
  EA^{[2n]}(B_\sigma^{[2n]}(s,\mathbf{z}))\leq 
  \max_{z\in\Delta^q}\{EA(\sigma(\mathbf{z}))\}+\tfrac{C(\sigma)}{2n}< c_2-\epsilon+\tfrac{C(\sigma)}{2n},
\end{equation*}
and\begin{equation*}
  EA^{[2n]}(B_\sigma^{[2n]}(1,\mathbf{z}))\leq 
  \max_{z\in\partial\Delta^q}\{EA(\sigma(\mathbf{z}))\}+\tfrac{C(\sigma)}{2n}< c_1-\epsilon+\tfrac{C(\sigma)}{2n}.
\end{equation*}
Notice that $C(\sigma)$ depends on $L$ and $\sigma$, we denote by $\bar{n}(L,\sigma):=\lceil\tfrac{C(\sigma)}{2\epsilon} \rceil$.
For each $n\geq \bar{n}(L,\sigma)$, $B_\sigma^{[2n]}$ satisfies the properties of the Lemma. 
$\hfill \qedsymbol$

~\\
 We also need the following lemma. 
See \cite[Lemma 1]{bangert_homology_1983}.
\begin{lem}\label{sec7:4}
  Let $(X,A)$ be a pair of topological space and $\alpha$ a singular relative $p$-cycle of $(X,A)$. Let $\Sigma(\alpha)$
  denotes the set of singular simplices of $\alpha$ together with all their faces. Suppose to every $\sigma\in\Sigma(\alpha)$,
  $\sigma:\Delta^q\rightarrow X,\ 0\leq q\leq p,$ there is assigned a map $P_\sigma: [0,1]\times\Delta^q\rightarrow X$
  such that 
  \begin{align*}
  (\textnormal{a})& P_\sigma(0,\cdot)=\sigma\text { for } z\in\Delta^q\\
(\textnormal{b})& P_\sigma(1,\Delta^q)\subset A\\
(\textnormal{c})& P_\sigma(s,\cdot)=\sigma \text{ if } \sigma(\Delta^q)\subset A\\
(\textnormal{d})& P_\sigma(\cdot,F_i(\cdot))=P_{\sigma\circ F_i}  \text{ for } i=0,\dots,q,
\text{ where } F_i:\Delta^{q-1}\rightarrow \Delta^q \text{ is the standard affine map  } \\
&\text{ onto the }i^{th} \text{ face of } \Delta^q.
  \end{align*}
  Then the homology class $[\alpha]\in H_p(X,A)$ vanishes.
\end{lem}

~\\
\noindent\emph{Proof of Lemma \ref{sec6:7}.}
By \eqref{sec6:6}, we have
$\rho_*^{}[\mu]\in H_*((EA_{T}^{})_{c+\epsilon},(EA_{T}^{})_c)$
for each $T>\max\{U,\tilde{T}(a,1)\}$. 
Besides, since $C_N(EA_{U,\mathcal{U}},\gamma)\neq 0$, $\gamma$ can not be a local minimum of $EA_{U,\mathcal{U}}$ on $\mathcal{X}$.
By the density of $\mathcal{X}$ in $EH(1)$, $(EA_{T})_c\neq \emptyset $ for each $T>\max\{U,\tilde{T}(a,1)\}$.

% Since $(EA^{}_{U,k})_c\neq\emptyset$, we have $(EA^{})_c\neq\emptyset.$
 Let $\Sigma(\mu)$ be the set of singular simplices of $\mu$ together with
all their faces. 
For each $\sigma:\Delta^q\rightarrow EH(1)$ in $\Sigma(\mu)$,
we notice that $\sigma(\Delta^q)\subset \mathcal{U}$, 
$\sigma$ has $C^{1}$-bounded image and
\begin{equation*}
  \sup_{z\in\Delta^q}\max_{t\in S_1}\{|\tfrac{d}{dt}\sigma(z)(t)|_{\sigma(z)(t)}\}<U,\quad \forall\sigma\in \Sigma(\mu).
\end{equation*}
Therefore $\rho^{}(\sigma)$ has $C^{1}$-bounded image with upper bound $U$. 
For each $T\geq U$, we have 
\begin{equation}\label{sec7:7}
  EA^{}(\rho^{}(\sigma)(z))=EA_{T}(\rho(\sigma)(z))<c+\epsilon,\quad \forall z\in \Delta^q,
\end{equation}
and thus 
\begin{equation*}
  \rho(\sigma)(\Delta^q)\subset (EA)_{c+\epsilon},\quad \rho(\sigma)(\partial \Delta^q)\subset (EA)_c.
\end{equation*}

Now we prove that
there exists $2\bar{n}=2\bar{n}(L_\theta,\sigma)\in\mathbb{K'}$
and a map
\begin{equation*}
  P_\sigma^{[2\bar{n}]}:[0,1]\times \Delta^q\rightarrow (EA^{[2\bar{n}]})_{c+\epsilon}
\end{equation*}
such that
\begin{align*}
  (\text{i})& P^{[2\bar{n}]}_\sigma(0,\cdot)=\psi^{[2\bar{n}]}\circ\rho^{}(\sigma)\\
(\text{ii})& P^{[2\bar{n}]}_\sigma(1,\Delta^q)\subset (EA^{[2\bar{n}]})_{c}\\
(\text{iii})& P^{[2\bar{n}]}_\sigma(s,\cdot)=\psi^{[2\bar{n}]}\circ\rho^{}(\sigma)\quad \text{ if } \rho(\sigma)(\Delta^q)\subset (EA)_{c}\\
(\text{iv})& P^{[2\bar{n}]}_\sigma(\cdot,F_i(\cdot))=P^{[2\bar{n}]}_{\sigma\circ F_i}\quad  \text{ for } i=0,\dots,q
%\text{ where } F_i:\Delta^{q-1}\rightarrow \Delta^q \text{ is the standard affine map  } \\&\text{ onto the }i^{th} \text{ face of } \Delta^q.
  \end{align*}
  It is noted that the image of $P^{[2\bar{n}]}_\sigma$ contained in $(EA^{[2\bar{n}]})_{c+\epsilon}$ implicitly requires the reversibility of
  $P^{[2\bar{n}]}_\sigma$ in time, i.e. 
  \begin{equation*}
    P^{[2\bar{n}]}_\sigma(\cdot,\cdot)(t)=P^{[2\bar{n}]}_\sigma(\cdot,\cdot)(-t),\quad \forall t\in \mathbb{R}.
  \end{equation*}
  In this case, the Bangert homotopy modified in Lemma \ref{sec7:3} is adequate for the construction of $P^{[2\bar{n}]}_\sigma$.
   Now we sketch the proof of the existence of $\bar{n}(L_\theta,\sigma)$ and $P^{[2\bar{n}]}_\sigma$ for each $\sigma\in\Sigma(\mu)$.
  We refer the reader to \cite{bangert_homology_1983,mazzucchelli_lagrangian_2011} for more details.

  If $\sigma\in\Sigma(\mu)$ is a trivial $q$-singular simplex, i.e. $\rho(\sigma)(\Delta^q)\subset (EA)_{c_1},$  
  we take $\bar{n}=1$ and $P^{[2]}_\sigma=\psi^{[2]}\circ\rho(\sigma).$ In the following, we always assume each $\sigma\in\Sigma(\mu)$ is 
  nontrivial. We prove the existence of $\bar{n}(L_\theta,\sigma)$ and $P_\sigma^{[2\bar{n}]}$ by induction on $p=\max\{\dim \sigma\,|\, \sigma\in\Sigma(\mu)\}$.
 
We only consider the case that $\sigma(\Delta^q)\nsubseteq (EA^{})_{c }$. 
If $p=0$ and $\eta:=\rho^{}(\sigma)(0)\in (EA^{})_{c+\epsilon}$ is a symmetric loop.
Since $(EA^{})_c\neq\emptyset$, we choose a path $\Gamma:[0,1]\rightarrow (EA^{})_{c+\epsilon}$ connecting $\eta$ and a loop in  $(EA^{})_c$.
We choose $\bar{n}=1$ and $P^{[2]}_\sigma(t,{0})=\psi^{[2]}\circ\Gamma(t)$. 

When $\mu$ is a relative $p$-cycle with $p>0$. For each $q$-singular simplex $\nu\in\Sigma(\mu)$ with $0\leq q<p$, we assume
$ P^{[2n]}_\nu$ is obtained for 
$2n\in\mathbb{K'}$ such that
\begin{equation*}
  2n\geq 2\tilde{n}(L,\mu):=\max\{\bar{n}(L,\nu)\in\mathbb{K'} \text{ for }q\text{-simplex }\nu \text{ with } 0\leq q<p\}.
\end{equation*}
Let $\sigma\in\Sigma(\mu)$ be a $p$-singular simplex, and $\sigma\circ F_i: \Delta^{p-1}\rightarrow (EA)_{c_2}$ be the $i^{th}$ face of $\sigma$.
Then $P^{[2n]}_{\sigma\circ F_i}$ is obtained for each face $\sigma\circ F_i$. 
We construct a new map 
  \begin{equation*}
    Q^{[2n]}: \{0\}\times \Delta^p\cup[0,\tfrac{1}{2}]\times \partial\Delta^p \rightarrow (EA^{[2n]})_{c+\epsilon}\label{sec7:2}
  \end{equation*}
  as \begin{equation*}
    Q^{[2n]}(s,z)=
    \begin{cases}
      \psi^{[2n]}\circ\rho(\sigma)(s,z)& \text{ if } (s,z)\in\{0\}\times\Delta^p,\\
      P^{[2n]}_{\sigma\circ F_i}(2s,z) & \text{ if } (s,F_i(z))\in [0,\tfrac{1}{2}]\times\partial\Delta^p,\ \forall i=0,\dots,p.
    \end{cases}
  \end{equation*}
  Since $\{0\}\times \Delta^p\cup [0,\tfrac{1}{2}]\times \partial\Delta^p$ is a retraction of $[0,\tfrac{1}{2}]\times \Delta^p$,
  $Q^{[2n]}$ can be extended as a continuous map 
  \begin{equation*}
    Q^{[2n]}:[0,\tfrac{1}{2}]\times\Delta^p\rightarrow (EA^{[2n]})_{c+\epsilon}.
  \end{equation*}
  By the property (ii) of $P^{[2n]}_{\sigma\circ F_i}$,
  we have 
  $Q^{[2n]}(\tfrac{1}{2},\partial\Delta^p)\subset (EA^{[2n]})_{c}.$ Therefore 
 \begin{equation*}
  \sigma':=Q^{[2n]}(\tfrac{1}{2},\cdot):\ (\Delta^p,\partial\Delta^p)\rightarrow ((EA^{[2n]})_{c+\epsilon},(EA^{[2n]})_{c})
 \end{equation*}
 is a $p$-singular simplex. By Lemma \ref{sec7:3} there exist
 $\bar{n}(L_\theta,\sigma')$ and a Bangert homotopy 
 \begin{equation*}
  B_{\sigma'}^{[2n]}: [0,1]\times (\Delta^p,\partial\Delta^p)\rightarrow ((EA^{[2n]})_{c+\epsilon},(EA^{[2n]})_{c}) 
 \end{equation*}
 for each $n\geq\bar{n}(L_\theta,\sigma')$ in $\mathbb{K'}$.

 Let $\bar{n}:=\bar{n}(L_\theta,\sigma)=\max\{\tilde{n}(L_\theta,\mu),\bar{n}(L_\theta,\sigma')\}$, then $P_\sigma^{[2\bar{n}]}$
  constructed by
  \begin{equation*}
  P_\sigma^{[2\bar{n}]}(s,z)=
  \begin{cases}
    Q^{[2\bar{n}]}(s,z)& \text{ if } (s,z)\in [0,\frac{1}{2}]\times\Delta^p,\\
    B^{[2\bar{n}]}_{\sigma'}(2s-1,z)&\text{ if } (s,z)\in [\tfrac{1}{2},1]\times\Delta^p.
  \end{cases}
\end{equation*}
satisfies the properties (i)(ii)(iii)(iv).

 Now we have obtained $\bar{n}(L_\theta,\sigma)$ and $P_\sigma^{[2\bar{n}]}$ for each $\sigma\in\Sigma(\mu)$. 
 Since each singular simplex $\sigma\in\Sigma(\mu)$ has $C^{1}$-bounded image, the Bangert homotopy $B_{\sigma}^{[2n]}$
  has $W^{1,\infty}$-bounded image.
  Therefore, the construction of $P_\sigma^{[2\bar{n}]}$ implies that $P_\sigma^{[2\bar{n}]}:[0,1]\times \Delta^q\rightarrow (EA^{[2\bar{n}]})_{c+\epsilon}$ has $W^{1,\infty}$-bounded image, i.e.
  \begin{equation*}
   (\text{v})\, \sup_{(t,z)\in [0,1]\times\Delta^q}\text{ess}\sup_{t\in S_1}\big|\tfrac{d}{dt}P^{[2\bar{n}]}_\sigma(s,z)(t)\big|_{P^{[2\bar{n}]}_\sigma(s,z)(t)}<T(L_{\theta},\sigma)
  \end{equation*}
  for some $T(L_{\theta},\sigma)>0.$ (i) implies that $T(L_{\theta},\sigma)>U$. Since $\Sigma(\mu)$ is a finite set, let 
  \begin{align*}
    \bar{m}&:=\bar{m}(L_\theta,[\mu])=\max\{2\bar{n}(L_{\theta},\sigma)\, |\, \sigma\in\Sigma(\mu)\},\\
    \bar{T}&:=\bar{T}(L_\theta,[\mu])=\max\{\tilde{T}(a,\bar{m}),\,T(L_{\theta},\sigma)\, |\, \sigma\in\Sigma(\mu)\}.
  \end{align*}
  Notice that for each $n\in\{n\in\mathbb{K'}\,|\, n\geq\bar{n}\}$, $P^{[2n]}_\sigma$ also satisfies the above properties (i)(ii)(iii)(iv). 
  Then for each $q$-singular simplex $\sigma\in\Sigma(\mu)$ and each
   $T\geq\bar{T}$, we have
  \begin{equation*}
    P_\sigma^{[\bar{m}]}([0,1]\times\Delta^q)\subset (EA_T^{[\bar{m}]})_{c+\epsilon}.
  \end{equation*}
  It is easy to verify that $P_\sigma^{[\bar{m}]}$ satisfies the properties
  \begin{align*}
    (\text{i'})& P^{[\bar{m}]}_\sigma(0,\cdot)=\psi^{[\bar{m}]}\circ\rho^{}(\sigma)\\
  (\text{ii'})& P^{[\bar{m}]}_\sigma(1,\Delta^q)\subset (EA_T^{[\bar{m}]})_{c}\\
  (\text{iii'})& P^{[\bar{m}]}_\sigma(s,\cdot)=\psi^{[\bar{m}]}\circ\rho^{}(\sigma)\quad \text{if } \rho(\sigma)(\Delta^q)\subset (EA_{T}^{})_{c}\\
  (\text{iv'})& P^{[\bar{m}]}_\sigma(\cdot,F_i(\cdot))=P^{[\bar{m}]}_{\sigma\circ F_i}\quad  \text{for } i=0,\dots,q
  %\text{ where } F_i:\Delta^{q-1}\rightarrow \Delta^q \text{ is the standard affine map  } \\&\text{ onto the }i^{th} \text{ face of } \Delta^q.
    \end{align*}
  Since $\psi_*^{[\bar{m}]}\circ\rho_*[\mu]\in H_*((EA^{[\bar{m}]}_{T})_{c+\epsilon},(EA^{[\bar{m}]}_{T})_c)$, by Lemma \ref{sec7:4}, we have 
  \begin{equation*}
    \psi_*^{[\bar{m}]}\circ\rho_*^{}[\mu]=0 \text{ in } H_*((EA_{T}^{[\bar{m}]})_{c+\epsilon},(EA_{T}^{[\bar{m}]})_c).
  \end{equation*}
Lemma \ref{sec6:7} is proved.

  \appendix   %仅一个附录时用appendix,否则\appendices
 % \setcounter{table}{0}   %从0开始编号,显示出来表会A1开始编号
 % \setcounter{figure}{0}
  %定义编号格式，在数字序号前加字符“A"
 % \renewcommand{\thetable}{A\arabic{table}}
   \renewcommand{\theequation}{A\arabic{equation}}
  \makeatletter
\@addtoreset{equation}{section}
\makeatother

\section{Appendix}
We introduce the generalized splitting lemma.
\begin{conditionsofSPL}
  Let $H$ be a Hilbert space with inner product $\langle\cdot,\cdot\rangle_H$ and the induced norm 
  $\Vert\cdot\Vert$. Let $X$ be a Banach space with norm $\Vert\cdot\Vert$, such that 

  (S) $X\subset H$ is dense in $H$ and the inclusion $X\hookrightarrow H$ is continuous, i.e. we may assume 
  $\Vert x\Vert\leq \Vert x\Vert_X, \, \forall x\in X$.
    
  We denote by $V^X:=V\cap X$ an open neighborhood of 0 in $X$, where $V$ is an open neighborhood of the origin 
  $0\in H$. 

  (F) Suppose that a functional $\mathcal{L}:V\rightarrow\mathbb{R}$ satisfies the following conditions:

  (F1) $\mathcal{L}:V\rightarrow \mathbb{R}$ is a $C^1$-functional. 
  The origin $0\in X$ is a critical point of $\mathcal{L}|_{V^X}$ and $\mathcal{L}$.

  (F2) There exists a map $\mathcal{A}\in C^1(V^X,X)$ such that
  \begin{equation*}
    d\mathcal{L}(x)(u)=\langle \mathcal{A}(x),u\rangle_H,\quad \forall x\in V^X,\, u\in X.
  \end{equation*} 

  (F3) There exists a map $\mathcal{B}\in C(V^X,L_s(H))$ such that 
  \begin{equation*}
    \langle d\mathcal{A}(x)[u],v\rangle_H=\langle \mathcal{B}_T(\gamma)u,v\rangle_H,\quad \forall x\in V^X,\, u,v\in X.
  \end{equation*}

(B) $\mathcal{B}$ has to satisfy the following properties:

(B1) 0 is either not in the spectrum $\sigma(\mathcal{B}(0))$ or is an isolated point of $\sigma(\mathcal{B}(0))$.

  (B2) If $u\in H$ such that $\mathcal{B}(0)(u)=v$ for some $v\in X$, then $u\in X$.
   
  (B3) The map $\mathcal{B}: V_X\rightarrow L_s(H,H)$ has a deformation 
  \begin{equation*}
    \mathcal{B}(x)=P(x)+Q(x),\quad x\in V\cap X
  \end{equation*}
  where $P(x):H\rightarrow H$ is a positive definite linear operator and $Q(x):H\rightarrow H$ is a compact 
  linear operator with the following properties

  (i) All eigenfunctions of the 
  operatr $\mathcal{B}(0)$ that correspond to negative eigenvalue belong to $X$;

  (ii) For any sequence $\{x_k\}\subset V\cap X$ with $\Vert x_k\Vert\rightarrow 0$ it holds that $\Vert P(x_k)u-P(0)u\Vert\rightarrow0$
  for any $u\in H$;

  (iii) The map $Q$ is continuous at $0$ with respect to the topology induced form $H$ on $V\cap X$;

  (iv) There exist positive constants $\epsilon>0$ and $C>0$ such that 
  \begin{equation*}
    \langle P(x)u,u\rangle\geq C\Vert u\Vert^2\quad \forall u\in H,\, \forall x\in B_H(0,\epsilon)\cap X.
  \end{equation*}
\end{conditionsofSPL}

\begin{rmk}
  The condition (B2)(B3) imply that both $H^0$ and $H^-$ are finitely dimensional subspaces contained in 
  $X$ (See \cite[Proposition B.2]{lu_splitting_2013}). Let $*=+,-,0.$ We denote by $P^*$ the orthogonal projections 
  onto $H^*,$ and by $X^*=X\cap H^*=P^*(H)$. Let $H^\pm=H^+\oplus H^-$ and $X^\pm:=X\cap H^\pm=X^-+H^+\cap X.$
  We denote by $B_\epsilon$ the neighboehood of $0\in H$ with radius $\epsilon$ and $B_\epsilon^*:=B_\epsilon\cap H^*,*=0,+,-$.
Note that $\dim H^-,\,\dim H^0 $ are the Morse index and nullity of $\mathcal{L}$ at critical point $0$, i.e. 
$m^-(\mathcal{L},0)=\dim H^-$ and $m^0(\mathcal{L},0)=\dim H^0$. 
\end{rmk}

\begin{thm}[{\cite[Theorem 1.1]{lu_corrigendum_2009}}]\label{secA:1}
  Let $\mathcal{L}:V\rightarrow\mathbb{R}$ be a functional satisfying the assumptions 
  (S)(F)(B) and $\nu>0$. There exists a positive $\epsilon\in\mathbb{R}$, a (unique) $C^1$ 
  map $h: B_{H^0}(0,\epsilon)\rightarrow X^\pm$ satisfying $h(0)=0$ and 
  \begin{equation*}
    (I-P^0)\mathcal{A}(z+h(z))=0\quad z\in B_\epsilon^0,
  \end{equation*}
  an open neighborhood $W$ of 0 in $H$ and an origin-preserving homeomorphism 
  \begin{equation*}
    \Phi: B_{H^0}(0,\epsilon)\times (B_\epsilon^++B_\epsilon^-)\rightarrow W
  \end{equation*}
  of form $\Phi(z,u^++u^-)=z+h(z)+\phi_z(u^++u^-)$ with $\phi_z(u^++u^-)\in H^\pm$ such that
  \begin{equation*}
    \mathcal{L}\circ\Phi(z,u^++u^-)=\Vert u^+\Vert^2-\Vert u^-\Vert^2+\mathcal{L}{(z+h(z))}
  \end{equation*}
  for all $(z,u^++u^-)\in B_\epsilon^0\times (B_\epsilon^++B_\epsilon^-),$ and that 
  \begin{equation*}
  \Phi(B_{H^0}(0,\epsilon))\times(B_\epsilon^++B_\epsilon^-)\subset X.
\end{equation*}
Moreover, the homeomorphism $\Phi$ also has properties:

(a) For each $z\in B_\epsilon^0$, $\Phi(z,0)=z+h(z),\phi_z(u^++u^-)\in H^-$ if and only if $u^+=0$.

(b) The restriction of $\Phi$ to $B_\epsilon^0\times B_\epsilon^-$ is a homeomorphism from 
$B_\epsilon^0\times B_\epsilon^-\subset X\times X$ onto $\Phi(B_\epsilon^0\times B_\epsilon^-)\subset X$
even if the topologies on these two sets are chosen as the indued one by $X$.

The map $h$ and the function $B_\epsilon^0\ni z\mapsto \beta(z):=\mathcal{L}(z+h(z))$ also satisfy:

(i) The map $h$ is $C^1$ in $V_\rho^X$ and 
\begin{equation*}
  h'(z)=-[(I-P^0)\mathcal{A'}(z)|_{X^\pm}]^{-1}(I-P^0)\mathcal{A'}(z),\quad \forall z\in B_\epsilon^0;
\end{equation*}

(ii) $ \beta$ is $C^{2}$ and
\begin{equation*}
  d\beta(z_0)(z)=\Angles{\mathcal{A}(z_0+h(z_0))}{z},\quad\forall z_0\in B_\epsilon^0,\, z\in H^0,
\end{equation*}
%and $d\beta$ is strictly F-differentiable at $0\in H^0$ and $d^2\beta(0)=0$;

(iii) If $\mathbf{0} $ is an isolated critical point of $\mathcal{L}|_{V^X}$, then it is also an isolated critical point of $\mathcal{L}^0$.
\end{thm}

\begin{cor}[Shifting, {\cite[Corollary 2.6]{lu_splitting_2013}}]\label{secA:3}
  Under the assumption of Theorem \ref{secA:1}, if 0 is an isolated critical point of $\mathcal{L}$, for each Abelian group $\mathbb{G}$ it holds that 
  \begin{equation*}
    C_*(\mathcal{L},\mathbf{0};\mathbb{G})=C_{*-m^-(\mathcal{L},0)}(\beta,\mathbf{0};\mathbb{G}).
  \end{equation*}
\end{cor}
\begin{thm}[{\cite[Corollary 2.11]{lu_splitting_2013}}]\label{secA:2}
  Under the assumption of Theorem \ref{secA:1}, let $\mathbf{0}$ be an isolated critical point of $\mathcal{L}$ with critical value $c:=\mathcal{L}(\mathbf{0})$, 
  then for any open neighborhood $W$ of $\mathbf{0}$ in $V$, the inclusion map 
  \begin{equation*}
    I_{X,H}:((\mathcal{L}_{X}|_{{W}^X})_c\cup\{\mathbf{0}\}, (\mathcal{L}_{X}|_{{W}^X})_c) \rightarrow (((\mathcal{L}|_{{W}})_c\cup\{\mathbf{0}\}), (\mathcal{L}|_{{W}})_c)
  \end{equation*}
  induces the following homology isomorphism with coefficients in a field $\mathbb{F}$
  \begin{equation*}
    I_{X,H*}: C_*(\mathcal{L}_{X},\mathbf{0};\mathbb{F})\xrightarrow{\simeq} C_*(\mathcal{L},\mathbf{0};\mathbb{F}).
  \end{equation*}
\end{thm}

 \bibliography{Brakeorbits}

\end{document}